\RequirePackage{ifpdf}
\ifpdf 
\documentclass[pdftex]{sigma}
\else
\documentclass{sigma}
\fi

\numberwithin{equation}{section}

\newtheorem{Theorem}{Theorem}[section]
\newtheorem{Corollary}[Theorem]{Corollary}
\newtheorem{Lemma}[Theorem]{Lemma}
\newtheorem{Proposition}[Theorem]{Proposition}
{\theoremstyle{definition}
\newtheorem{Note}[Theorem]{Note}
\newtheorem{Definition}[Theorem]{Definition}
\newtheorem{Problem}[Theorem]{Problem}
}

\begin{document}

\allowdisplaybreaks

\renewcommand{\thefootnote}{$\star$}

\renewcommand{\PaperNumber}{069}

\FirstPageHeading

\ShortArticleName{The Universal Askey--Wilson Algebra}

\ArticleName{The Universal Askey--Wilson Algebra\footnote{This paper is a
contribution to the Special Issue ``Relationship of Orthogonal Polynomials and Special Functions with Quantum Groups and Integrable Systems''. The
full collection is available at
\href{http://www.emis.de/journals/SIGMA/OPSF.html}{http://www.emis.de/journals/SIGMA/OPSF.html}}}

\Author{Paul TERWILLIGER}

\AuthorNameForHeading{P.~Terwilliger}

\Address{Department of Mathematics, University of Wisconsin,  Madison, WI 53706-1388, USA}
\Email{\href{mailto:terwilli@math.wisc.edu}{terwilli@math.wisc.edu}}

\ArticleDates{Received April 17, 2011, in f\/inal form July 09, 2011;  Published online July 15, 2011}

\Abstract{Let $\mathbb F$ denote a f\/ield, and f\/ix a nonzero $q \in \mathbb F$ such
that $q^4\not=1$.
We def\/ine an associative $\mathbb F$-algebra $\Delta=\Delta_q$ by generators and
relations in the following way.
The generators are
$A$, $B$, $C$. The relations assert that each of
\begin{gather*}
A + \frac{qBC-q^{-1}CB}{q^2-q^{-2}},
\qquad
B + \frac{qCA-q^{-1}AC}{q^2-q^{-2}},
\qquad
C + \frac{qAB-q^{-1}BA}{q^2-q^{-2}}
\end{gather*}
is central in $\Delta$.
We call $\Delta$ the  {\it universal Askey--Wilson algebra}.
We discuss how $\Delta$ is related to the original
Askey--Wilson
algebra AW(3) introduced by A.~Zhedanov.
Multiply  each of the above central elements by
 $q+q^{-1}$ to obtain
  $\alpha$,
  $\beta$,
  $\gamma$.
We give an alternate presentation for $\Delta$
by generators and relations;
the generators are $A$, $B$, $\gamma$.
We give a faithful action of the
modular group
${\rm  {PSL}}_2(\mathbb Z)$  on $\Delta$ as a group of automorphisms;
one generator sends $(A,B,C)\mapsto (B,C,A)$
and another generator sends $(A,B,\gamma) \mapsto (B,A,\gamma)$.
We show that
$\lbrace A^iB^jC^k \alpha^r\beta^s\gamma^t| i,j,k,r,s,t\geq 0\rbrace$
is a basis for the
$\mathbb F$-vector space $\Delta$.
We show that the center $Z(\Delta)$ contains the
element
\begin{gather*}
\Omega =
qABC+q^2A^2+q^{-2}B^2+q^2C^2-qA\alpha-q^{-1}B\beta -q C\gamma.
\end{gather*}
Under the assumption that $q$ is not a root of unity,
we
show that $Z(\Delta)$ is generated by~$\Omega$, $\alpha$, $\beta$, $\gamma$
and that $Z(\Delta)$ is isomorphic to a polynomial algebra in 4
variables.
Using the alternate presentation we  relate $\Delta$ to
the $q$-Onsager algebra.
We describe the 2-sided ideal
$\Delta\lbrack \Delta,\Delta\rbrack \Delta$ from several points of view.
Our main result here is that
$\Delta\lbrack \Delta,\Delta \rbrack \Delta + \mathbb F 1$ is equal to the intersection
of $(i)$~the subalgebra of $\Delta$ generated by $A$, $B$;
$(ii)$~the subalgebra of $\Delta$ generated by $B$, $C$;
$(iii)$~the subalgebra of $\Delta $ generated by $C$, $A$.}


\Classification{33D80; 33D45}

\renewcommand{\thefootnote}{\arabic{footnote}}
\setcounter{footnote}{0}

\section{Introduction}

In \cite{Zhidd}
A.~Zhedanov introduced the Askey--Wilson algebra
${\rm AW}=
{\rm AW}(3)$
and used it to
describe the Askey--Wilson polynomials~\cite{Askey}.
Since then, ${\rm AW}$ has become one of the main
objects in the theory
of the Askey scheme
of orthogonal polynomials~\cite{GYLZmut,
  GYZlinear,
  GH1,
KoeSwa, Koo1,Koo2, VZ2,VZ3}.
It is particularly useful
in the theory of
Leonard pairs~\cite{NT:span,
LS99,
madrid,
aw,
vidunas,
V2}
and Leonard triples~\cite{Curt2,Curt1,
koro}.
The algebra
${\rm AW}$
is related to
the algebra
$U_q(\mathfrak{sl}_2)$~\cite{GYZnature,
  GYZlinear,
Rosengren,
rosengren,
wieg}
and the algebra
{$U\sb q(\rm su\sb 2)$}
\cite{atak1, atak2,atak3}.
There is a connection to the
double af\/f\/ine Hecke algebra
of type $(C^\vee_1,C_1)$~\cite{dahater,
Koo1,Koo2}.
The $\mathbb Z_3$-symmetric quantum
algebra $O'_q({\rm so}_3)$~\cite[Remark~6.11]{ciccoli},
\cite{fairlie},
\cite[Section~3]{su2},
\cite{posta,
posta2,
jordan,
odesskii}
is a special case of
${\rm AW}$,  and
the recently introduced Calabi--Yau algebras~\cite{calabi} give a generalization
of
${\rm AW}$.
The algebra ${\rm AW}$
plays a role in
integrable systems
\cite{aneva,bas2,bas1,bas3,bas4,bas5,bas6,bas7,bas8,basnc,
lav1,lav2,lav3,VZ}
and quantum mechanics~\cite{odake2,
odake}, as well as
the theory of quadratic algebras~\cite{jordan,
kalnins,
ros}.
There is a classical version of
${\rm AW}$
that has
a Poisson algebra structure
\cite{GYLZmut}, \cite[equation~(2.9)]{koro}, \cite[equations (26)--(28)]{oblomkov}, \cite{LPcm}.

In this paper we introduce a
central extension of ${\rm AW}$
called the {\it universal Askey--Wilson algebra}.
This central extension, which we denote by
$\Delta$, is related to
${\rm AW}$ in the following way.
There is a
reduced $\mathbb Z_3$-symmetric presentation
of ${\rm AW}$ that involves three scalar
parameters besides $q$
\cite[equation~(6.1)]{koro}.
Up to normalization, the algebra
$\Delta $ is what one gets from
this presentation
by reinterpreting the three scalar parameters
as central elements in the algebra.
By construction
$\Delta$ has no scalar parameters besides $q$,
and there exists a surjective algebra homomorphism
$\Delta \to
{\rm AW}$.
One advantage of
$\Delta$ over ${\rm AW}$
is that $\Delta$ has a larger automorphism
group. Our def\/inition of $\Delta$ was inspired by
\cite[Section 3]{dahater}, which in turn was motivated by
\cite{ion}.

Let us now bring in more detail, and
recall the def\/inition of
${\rm AW}$. There are at least three presentations in the
literature; the original one involving three
generators
\cite[equations~(1.1a)--(1.1c)]{Zhidd},
one involving two generators
\cite[equations~(2.1), (2.2)]{Koo1},
\cite[Theorem~1.5]{aw},
and a $\mathbb Z_3$-symmetric presentation involving three generators
\cite[equation~(6.1)]{koro},
\cite[p.~101]{ros},
\cite{smithandbell},
\cite[Section~4.3]{wieg}.
We will use the presentation
in \cite[equation~(6.1)]{koro},
although we adjust the
normalization and replace
$q$ by $q^2$
in order to illuminate the underlying symmetry.

Our conventions for the paper are as follows.
An algebra is meant to be associative and have a 1.
A subalgebra has the same 1 as the parent algebra.
We f\/ix a f\/ield $\mathbb F$
and  a nonzero $q \in \mathbb F$ such that $q^4\not=1$.

\begin{Definition}[\protect{\cite[equation~(6.1)]{koro}}]
\label{def:aw3}
Let $a$, $b$, $c$ denote scalars in $\mathbb F$.
Def\/ine the $\mathbb F$-algebra ${\rm AW}={\rm AW}_q(a,b,c)$
 by generators and relations in the following way.
The generators are~$A$,~$B$,~$C$.
The relations assert that
\begin{gather*}
 A+ \frac{qBC-q^{-1}CB}{q^2-q^{-2}},
\qquad B+
\frac{qCA-q^{-1}AC}{q^2-q^{-2}},
\qquad
C+
\frac{qAB-q^{-1}BA}{q^2-q^{-2}}
\end{gather*}
are equal to
$
a/(q+q^{-1})$,
$b/(q+q^{-1})$,
$c /(q+q^{-1})$
respectively. We call ${\rm AW}$ the
{\it Askey--Wilson algebra} that corresponds to
$a$, $b$, $c$.
\end{Definition}

 We now introduce the algebra $\Delta$.

\begin{Definition}
\label{def:uaw}
Def\/ine an $\mathbb F$-algebra
$\Delta=\Delta_q$
 by generators and relations in the following way.
The generators are~$A$,~$B$,~$C$.
The relations assert that each of
\begin{gather}
 A+ \frac{qBC-q^{-1}CB}{q^2-q^{-2}},
\qquad B+
\frac{qCA-q^{-1}AC}{q^2-q^{-2}},
\qquad
C+
\frac{qAB-q^{-1}BA}{q^2-q^{-2}}
\label{eq:comlist}
\end{gather}
is central in $\Delta$.
We call $\Delta$ the {\it universal Askey--Wilson algebra}.
\end{Definition}

\begin{Definition}
\label{def:abc}
\rm
For the
three central elements in
(\ref{eq:comlist}), multiply each by $q+q^{-1}$ to get~$\alpha$,~$\beta$,~$\gamma$. Thus
\begin{gather}
A+ \frac{qBC-q^{-1}CB}{q^2-q^{-2}}
=
\frac{\alpha}{q+q^{-1}},
\label{eq:u1}
\\
B+
\frac{qCA-q^{-1}AC}{q^2-q^{-2}}
=
\frac{\beta}{q+q^{-1}},
\label{eq:u2}
\\
C+
\frac{qAB-q^{-1}BA}{q^2-q^{-2}}
=
\frac{\gamma}{q+q^{-1}}.
\label{eq:u3}
\end{gather}
Note that each of $\alpha$, $\beta$, $\gamma$ is central in
$\Delta$.
(The purpose of the factor $q+q^{-1}$ is
to make the upcoming formula~(\ref{eq:caspreview})
more attractive.)
\end{Definition}

From the construction we obtain the following
result.

\begin{Lemma}
Let $a$, $b$, $c$ denote scalars in $\mathbb F$
and let
${\rm AW}$ denote the corresponding
Askey--Wilson algebra.
Then there exists a surjective $\mathbb F$-algebra
homomorphism $\Delta \to {\rm AW}$ that sends
\begin{gather*}
A \mapsto A,
\qquad
B \mapsto B,
\qquad
C \mapsto C,
\qquad
\alpha \mapsto a,
\qquad
\beta \mapsto b,
\qquad
\gamma \mapsto c.
\end{gather*}
\end{Lemma}

In this paper we begin a
comprehensive study of the algebra~$\Delta$. For now
we consider
the ring-theoretic aspects, and leave the
representation theory for some future paper.
Our main results are summarized as follows.
We give an alternate presentation for~$\Delta$
by generators and relations;
the generators are $A$, $B$, $\gamma$.
Following
\cite[Lemma~5.2]{Curt2},
\cite[Theorem~5.1]{dahater},
\cite{koro},
\cite[Section~1.2]{oblomkov}
we give a faithful action of the
modular group
${\rm  {PSL}}_2(\mathbb Z)$  on $\Delta$ as a group of automorphisms;
one generator sends $(A,B,C)\mapsto (B,C,A)$
and another generator sends $(A,B,\gamma) \mapsto (B,A,\gamma)$.
Following
\cite[Theorem~1]{posta2},
\cite[Proposition~6.6$(i)$]{jordan}
we show that
\begin{gather*}
A^iB^jC^k \alpha^r\beta^s\gamma^t,
\qquad   i,j,k,r,s,t \geq 0
\end{gather*}
is a basis for the $\mathbb F$-vector space
$\Delta$.
Following
\cite[Lemma~1]{posta2},
\cite[Proposition~3]{koro},
\cite[equa\-tion~(1.3)]{Zhidd}
we show that the center $Z(\Delta)$ contains a
Casimir element
\begin{gather}
\label{eq:caspreview}
\Omega =
qABC+q^2A^2+q^{-2}B^2+q^2C^2-qA\alpha-q^{-1}B\beta -q C\gamma.
\end{gather}
Under the assumption that $q$ is not a root of unity,
we
show that $Z(\Delta)$ is generated by
$\Omega$, $\alpha$, $\beta$, $\gamma$
and that $Z(\Delta)$ is isomorphic to a polynomial algebra in 4~variables.
Using the alternate presentation we  relate $\Delta$ to
the $q$-Onsager algebra
\cite[Section 4]{basnc},
\cite{IT:aug},
\cite[Def\/inition~3.9]{qSerre}.
We describe the 2-sided ideal
$\Delta\lbrack \Delta,\Delta\rbrack \Delta$ from several points of view.
Our main result here is that
$\Delta\lbrack \Delta,\Delta \rbrack \Delta + \mathbb F 1$
is equal to the intersection
of $(i)$~the subalgebra of $\Delta$ generated by~$A$,~$B$;
$(ii)$~the subalgebra of $\Delta$ generated by~$B$,~$C$;
$(iii)$~the subalgebra of $\Delta $ generated by~$C$,~$A$.
At the end of the paper we list some open problems
that are intended to
motivate further research.

\section[Another presentation of $\Delta$]{Another presentation of $\boldsymbol{\Delta}$}

A bit later in the paper
we will discuss automorphisms of~$\Delta $.
To facilitate this discussion
we give
another presentation for~$\Delta $ by generators and relations.

\begin{Lemma}
\label{lem:abgpres}
The $\mathbb F$-algebra $\Delta $ is generated by
$A$, $B$, $\gamma$. Moreover
\begin{gather}
C  =  \frac{\gamma}{q+q^{-1}}
-\frac{qAB-q^{-1}BA}{q^2-q^{-2}},
\label{eq:getc}
\\
\alpha  =  \frac{B^2A-(q^2+q^{-2})BAB+AB^2+
(q^2-q^{-2})^2A+(q-q^{-1})^2B\gamma}{(q-q^{-1})(q^2-q^{-2})},
\label{eq:getalpha}
\\
\beta  =  \frac{A^2B-(q^2+q^{-2})ABA+BA^2+
(q^2-q^{-2})^2B+(q-q^{-1})^2A\gamma}{(q-q^{-1})(q^2-q^{-2})}.
\label{eq:getbeta}
\end{gather}
\end{Lemma}

\begin{proof}\sloppy
Line~(\ref{eq:getc})
is from (\ref{eq:u3}).
To get
(\ref{eq:getalpha}),
(\ref{eq:getbeta})
eliminate $C$ in
(\ref{eq:u1}),
(\ref{eq:u2})
using
(\ref{eq:getc}).
\end{proof}

 Recall the notation
\begin{gather*}
\lbrack n \rbrack_q= \frac{q^n-q^{-n}}{q-q^{-1}}, \qquad
n=0,1,2,\ldots.
\end{gather*}

{\samepage

\begin{Theorem}
\label{prop:altpres}
The $\mathbb F$-algebra $\Delta $ has a  presentation by generators
$A$, $B$, $\gamma$ and
relations
 \begin{gather*}
 A^3B-\lbrack 3\rbrack_q A^2BA+\lbrack 3\rbrack_q ABA^2-BA^3
= -\big(q^2-q^{-2}\big)^2(AB-BA),
\\
 B^3A-\lbrack 3\rbrack_q B^2AB+\lbrack 3\rbrack_q BAB^2-AB^3
= -\big(q^2-q^{-2}\big)^2(BA-AB),
\\
 A^2B^2-B^2A^2+\big(q^2+q^{-2}\big)(BABA-ABAB) = -\big(q-q^{-1}\big)^2(AB-BA)\gamma,
\\
  \qquad \gamma A = A \gamma,   \qquad \gamma B=B \gamma.
 \end{gather*}
\end{Theorem}

\begin{proof}
Use Lemma~\ref{lem:abgpres} to
express the def\/ining relations
for $\Delta $ in terms of $A$, $B$, $\gamma$.
\end{proof}}

\begin{Note}
The f\/irst two equations
in Theorem~\ref{prop:altpres}
are known as the {\it tridiagonal relations}
\cite[De\-f\/i\-nition~3.9]{qSerre}.
These
relations
have appeared in algebraic combinatorics
\cite[Lemma 5.4]{tersub3},
 the theory of tridiagonal pairs
\cite{TD00,
IT:qRacah,
LS99,
qSerre,
madrid},
and integrable systems
\cite{bas2,bas1,bas3,bas4,bas5,bas6,bas7,bas8,basnc}.
\end{Note}

\section[An action of ${\rm  {PSL}}_2(\mathbb Z)$ on $\Delta$]{An action of $\boldsymbol{{\rm  {PSL}}_2(\mathbb Z)}$ on $\boldsymbol{\Delta}$}

We now consider
some automorphisms of~$\Delta $.
Recall that the modular group
${\rm  {PSL}}_2(\mathbb Z)$ has a~presentation
by generators~$\rho$, $\sigma$ and
relations $\rho^3=1$, $\sigma^2=1$.
See for example \cite{alpern}.

\begin{Theorem}
\label{thm:try}
The group
${\rm  {PSL}}_2(\mathbb Z)$ acts on
$\Delta $ as a group of automorphisms in the following way:
\[
\begin{array}{@{}c| ccc | c c c}
u &  A & B  & C
&  \alpha & \beta  & \gamma
\\
\hline
\rho(u) &  B & C & A
&  \beta  & \gamma & \alpha
\\
\sigma(u) &  B & A & C+\frac{AB-BA}{q-q^{-1}}
&  \beta & \alpha & \gamma
\end{array}
\]
\end{Theorem}

\begin{proof}
By Def\/inition
\ref{def:uaw} there exists an automorphism
$P$ of $\Delta $ that sends
\begin{gather*}
A \mapsto B,
\qquad
B \mapsto C,
\qquad
C \mapsto A.
\end{gather*}
Observe $P^3=1$. By
(\ref{eq:u1})--(\ref{eq:u3})
the map $P$ sends
\begin{gather*}
\alpha \mapsto \beta,
\qquad
\beta \mapsto \gamma,
\qquad
\gamma \mapsto \alpha.
\end{gather*}
By Theorem
\ref{prop:altpres} there exists
an automorphism
$S$ of $\Delta $ that sends
\begin{gather*}
A \mapsto B,
 \qquad
B \mapsto A,
\qquad
\gamma \mapsto \gamma.
\end{gather*}
Observe $S^2=1$. By
Lemma
\ref{lem:abgpres}
the map $S$
sends
\begin{gather*}
\alpha \mapsto \beta,
\qquad
\beta \mapsto \alpha,
\qquad
C \mapsto C + \frac{AB-BA}{q-q^{-1}}.
\end{gather*}
The result follows.
\end{proof}

In Theorem
\ref{thm:try}
we gave an
 action of ${\rm  {PSL}}_2(\mathbb Z)$  on $\Delta $.
Our next goal is to show that this
action is faithful.

Let $\lambda $ denote an indeterminate.
Let $\mathbb F\lbrack \lambda,\lambda^{-1}\rbrack$
denote the $\mathbb F$-algebra consisting of the Laurent
polynomials in $\lambda$ that have all coef\/f\/icients
in~$\mathbb F$.
We will be discussing the $\mathbb F$-algebra
\begin{gather*}
\Lambda ={\rm {Mat}}_2(\mathbb F) \otimes_{\mathbb F} \mathbb F\big\lbrack \lambda, \lambda^{-1}\big\rbrack.
\end{gather*}
We view elements of $\Lambda $ as $2\times 2$ matrices that have entries
in
$\mathbb F\lbrack \lambda, \lambda^{-1}\rbrack $.
From this point of view the product operation
for $\Lambda $ is ordinary matrix multiplication,
and the multiplicative identity in~$\Lambda $
is the identity matrix~$I$.
For notational convenience def\/ine
$\mu = \lambda + \lambda^{-1}$.

For later use we now describe the center $Z(\Lambda)$.{\samepage
\begin{Lemma}
\label{lem:zh}
For all $\eta \in \Lambda $ the following
$(i)$, $(ii)$
are equivalent:
\begin{enumerate}\itemsep=0pt
\item[$(i)$] $\eta  \in Z(\Lambda)$.
\item[$(ii)$] There exists $\theta \in
 \mathbb F\lbrack \lambda, \lambda^{-1} \rbrack $ such that $\eta=\theta I$.
\end{enumerate}
\end{Lemma}

\begin{proof}
Routine.
\end{proof}}

\begin{Definition}
\label{def:xyz}
Let $ \mathcal A$, $\mathcal B$, $\mathcal C$ denote the following elements of
$\Lambda $:
\begin{gather*}
 \mathcal A =
\left(
\begin{array}{  c c}
\lambda & 1-\lambda^{-1}   \\
0 & \lambda^{-1}
\end{array}
\right),
\qquad
\mathcal B =
\left(
\begin{array}{  c c}
\lambda^{-1} &  0 \\
\lambda-1 & \lambda
\end{array}
\right),
 \qquad \mathcal C =
\left(
\begin{array}{  c c}
1 &  1-\lambda \\
\lambda^{-1}-1 & \lambda+ \lambda^{-1}-1
\end{array}
\right).
\end{gather*}
\end{Definition}

\begin{Lemma}
\label{lem:xyz}
We have
\begin{gather*}
\mathcal A \mathcal B\mathcal C=I, \qquad
\mathcal A+{\mathcal A}^{-1}=\mu I,
\qquad
\mathcal B+{\mathcal B}^{-1}=\mu I,
\qquad
\mathcal C+{\mathcal C}^{-1}=\mu I.
\end{gather*}
\end{Lemma}

\begin{proof}
Use Def\/inition
\ref{def:xyz}.
\end{proof}

\begin{Lemma}
\label{lem:multtable2}
In the algebra $\Lambda $ the elements $\mathcal A$, $\mathcal B$, $\mathcal C$
multiply as follows:
\[
\begin{array}{@{}c| c c c}
 &  \mathcal A & \mathcal B  & \mathcal C \\
\hline
\mathcal A &   \mu \mathcal A  - I &   \mu I - \mathcal C
&  \mu \mathcal A  + \mathcal B + \mu \mathcal C  -\mu^2 I \\
\mathcal B
&
\mu \mathcal B + \mathcal C + \mu \mathcal A  -\mu^2 I
&
  \mu \mathcal B - I
&
 \mu I- \mathcal A
\\
\mathcal C
&
 \mu I - \mathcal B
&
\mu \mathcal C  + \mathcal A + \mu \mathcal B  -\mu^2 I
&
\mu \mathcal C - I
\end{array}
\]
\end{Lemma}

\begin{proof}
Use Lemma
\ref{lem:xyz}.
\end{proof}

 The algebra $\Lambda $ is not generated by
$\mathcal A$, $\mathcal B$, $\mathcal C$. However
we do have the following result.

\begin{Lemma}
\label{lem:whatittakes}
Suppose  $\eta  \in \Lambda $ commutes with at
 least two of
$\mathcal A$, $\mathcal B$, $\mathcal C$.
Then $\eta  \in Z(\Lambda )$.
\end{Lemma}

\begin{proof}
For $1 \leq i,j\leq 2$
let $\eta_{ij}$ denote the $(i,j)$-entry of $\eta$.
The matrix $\eta $ commutes with each of~$\mathcal A$,~$\mathcal B$,~$\mathcal C$
since $\mathcal A\mathcal B\mathcal C=I$.
In the equation $\eta \mathcal A=\mathcal A \eta  $,
evaluate $\mathcal A$ using Def\/inition~\ref{def:xyz}, and simplify the result to get
$\eta_{21}=0$. Similarly using $\eta \mathcal B=\mathcal B\eta $
we f\/ind $\eta_{12}=0$ and
$\eta_{11}=\eta_{22}$.
Therefore $\eta =\eta_{11} I\in Z(\Lambda )$.
\end{proof}

  Next we describe an
action of ${\rm  {PSL}}_2(\mathbb Z)$ on $\Lambda $ as
a group of automorphisms.

\begin{Definition}
\label{def:ps}
Let $p$ and $s$ denote the following elements
of $\Lambda $:
\begin{gather*}
p =
\left(
\begin{array}{  c c}
0 & -1   \\
1 & 1
\end{array}
\right),
\qquad
s =
\left(
\begin{array}{  c c}
0 &  1 \\
\lambda  & 0
\end{array}
\right).
\end{gather*}
\end{Definition}

\begin{Lemma}
\label{lem:4}
The following $(i)$--$(iv)$ hold.
\begin{enumerate}\itemsep=0pt
\item[$(i)$] $\det (p)=1$ and  $\det (s)=-\lambda$.
\item[$(ii)$]
 $p^3=-I$ and $s^2=\lambda I $.
\item[$(iii)$]
$p \mathcal A p^{-1} = \mathcal B$,
$p  \mathcal B p^{-1} = \mathcal C$,
$p \mathcal C p^{-1} = \mathcal A$.
\item[$(iv)$]
$s \mathcal A s^{-1} = \mathcal B$ and
$s \mathcal B s^{-1} = \mathcal A$.
\end{enumerate}
\end{Lemma}

\begin{proof}
$(i)$, $(ii)$
Use Def\/inition~\ref{def:ps}.
$(iii)$, $(iv)$ Use
Def\/inition~\ref{def:xyz}
and
Def\/inition~\ref{def:ps}.
\end{proof}

\begin{Lemma}
\label{lem:conj}
The group
${\rm  {PSL}}_2(\mathbb Z)$ acts on $\Lambda $ as a group of
automorphisms such
that
$\rho(\eta )= p \eta p^{-1}$
and
$\sigma (\eta ) = s \eta s^{-1}$ for all
$\eta  \in \Lambda $.
\end{Lemma}

\begin{proof}
By Lemma~\ref{lem:4}$(ii)$ the elements  $p^3$, $s^2$ are in $Z(\Lambda)$.
\end{proof}

\begin{Lemma}
\label{lem:faith1}
The action of
${\rm  {PSL}}_2(\mathbb Z)$ on $\Lambda $ is faithful.
\end{Lemma}

\begin{proof}
Pick an integer $n\geq 1$.
Consider an element $\eta \in \Lambda $
of the form
$\eta =\eta_1 \eta_2 \cdots \eta_n$ such that
 for $1 \leq i \leq n$,  $\eta_i=s$
for one parity of $i$ and
$\eta_i \in \lbrace p,p^{-1} \rbrace$ for the other parity of~$i$.
We show that $\eta \not\in Z(\Lambda)$.
To this end we assume
 $\eta \in Z(\Lambda )$ and get a contradiction.
Let $\ell$ denote the number of times~$s$ occurs among
$\lbrace \eta_i\rbrace_{i=1}^n$.
Assume for the moment $\ell = 0$.
Then $n=1$ so
$\eta =\eta_1 \in \lbrace p, p^{-1}\rbrace$.
The elements $p$, $p^{-1}$ are not in~$Z(\Lambda)$, for a contradiction.
Therefore $\ell \neq 0$.
From the nature of the matrices~$p$,~$s$ in Def\/inition~\ref{def:ps},
we may
view $\eta $ as a polynomial in $\lambda$
that has
coef\/f\/icients in
${\rm {Mat}}_2(\mathbb F)$ and degree at most $\ell$.
Call this polynomial~$f$.
We claim that the degree of~$f$ is exactly~$\ell$.
To prove the claim,
write
\begin{gather*}
s = s_0+ s_1 \lambda,   \qquad
s_0= \left(
\begin{array}{  c c}
0 & 1  \\
0 & 0
\end{array}
\right),
\qquad
s_1= \left(
\begin{array}{  c c}
0 & 0  \\
1 & 0
\end{array}
\right).
\end{gather*}
Let $m
\in {\rm {Mat}}_2(\mathbb F)$ denote the coef\/f\/icient of
$\lambda^\ell$ in $f$.
The matrix $m$ is obtained from
$\eta_1 \eta_2\cdots \eta_n$ by
replacing each occurrence of $s$ by $s_1$.
Using $s_1p s_1=-s_1$ and
 $s_1p^{-1} s_1=s_1$  we f\/ind
$m\in \lbrace \pm p^i s_1 p^j |
-1 \leq i,j\leq 1\rbrace $. The matrix $p$ is invertible
and $s_1\not=0$ so
$p^i s_1 p^j\not=0$ for $-1 \leq i,j \leq 1$.
Therefore $m\not=0$ and the claim is proved.
Let $\kappa $ denote the $(1,1)$-entry of the matrix $\eta $.
Then $\eta=\kappa I$ since $\eta \in Z(\Lambda )$.
In the equation
$\eta_1 \eta_2\cdots \eta_n=\kappa I$ take the determinant of
each side and use
Lemma~\ref{lem:4}$(i)$
to get $(-\lambda)^\ell = \kappa^2$.
Therefore $\ell $ is even and
$\kappa=\pm \lambda^{\ell/2}$.
Now $\eta =\pm \lambda^{\ell/2} I$, so
the above polynomial $f$ has degree $\ell/2$.
But $\ell \not=0$ so $\ell > \ell/2$
for a contradiction.
the result follows.
\end{proof}

We now display an algebra homomorphism $\Delta \to \Lambda $.

\begin{Lemma}
\label{prop:uhomoh}
There exists a unique $\mathbb F$-algebra homomorphism
$\pi:\Delta \to \Lambda $
 that sends
\begin{gather*}
A \mapsto q \mathcal A+q^{-1}{\mathcal A}^{-1},
 \qquad
B \mapsto q \mathcal B+q^{-1}{\mathcal B}^{-1},
\qquad
C \mapsto q\mathcal C+q^{-1}{\mathcal C}^{-1}.
\end{gather*}
The homomorphism $\pi$ sends
\begin{gather}
 \alpha \mapsto \nu I,
 \qquad
\beta \mapsto  \nu I,
 \qquad
\gamma \mapsto \nu I,
\label{eq:pisend}
\end{gather}
where $\nu= (q^2+q^{-2})\mu + \mu^2$.
\end{Lemma}

\begin{proof}
Def\/ine
\begin{gather*}
A^\vee = q\mathcal A+q^{-1}{\mathcal A}^{-1},
\qquad
B^\vee =q\mathcal B+q^{-1}{\mathcal B}^{-1},
\qquad
C^\vee = q\mathcal C+q^{-1}{\mathcal C}^{-1}.
\end{gather*}
By Lemma~\ref{lem:xyz}
and Lemma~\ref{lem:multtable2},
\begin{gather}
(q+q^{-1})A^\vee+ \frac{qB^\vee C^\vee-q^{-1}C^\vee B^\vee}{q-q^{-1}}
  = \nu I,
\label{eq:vee1}
\\
(q+q^{-1})B^\vee+
\frac{qC^\vee A^\vee -q^{-1}A^\vee C^\vee }{q-q^{-1}}
 =  \nu I,
\label{eq:vee2}
\\
(q+q^{-1})C^\vee +
\frac{qA^\vee B^\vee -q^{-1}B^\vee A^\vee }{q-q^{-1}}
 =  \nu I.
\label{eq:vee3}
\end{gather}
By
(\ref{eq:vee1})--(\ref{eq:vee3}) and
since
 $\nu I$ is central, the elements
$A^\vee$, $B^\vee$, $C^\vee$ satisfy the
def\/ining relations for~$\Delta $ from
Def\/inition~\ref{def:uaw}.
Therefore
the homomorphism $\pi$ exists. The homomorphism $\pi$
is unique since $A$, $B$, $C$ generate~$\Delta $. Line
(\ref{eq:pisend})
follows from
Def\/inition~\ref{def:abc} and
(\ref{eq:vee1})--(\ref{eq:vee3}).
\end{proof}

\begin{Lemma}
\label{lem:comdiagf}
For $g \in
{\rm  {PSL}}_2(\mathbb Z)$
 the following diagram commutes:
\[
\begin{CD}
\Delta  @>\pi>>
         \Lambda  \\
	  @Vg VV                     @VVg V \\
	 \Delta
	  @>\pi>>
	          \Lambda
		    \end{CD}
		    \]

\end{Lemma}

\begin{proof}
The elements
$\rho$, $\sigma$ form a generating set for
${\rm  {PSL}}_2(\mathbb Z)$; without loss
we may assume that~$g$ is contained in this set.
By Theorem~\ref{thm:try}
the action of~$\rho$ on~$\Delta $
cyclically permutes~$A$,~$B$,~$C$.
By Lemma~\ref{lem:4}$(iii)$
the action of $\rho$ on $\Lambda $
cyclically permutes
$\mathcal A$, $\mathcal B$, $\mathcal C$.
By Theorem~\ref{thm:try}
the action of~$\sigma$ on~$\Delta $ swaps $A$, $B$ and f\/ixes
$\gamma$. By Lemma~\ref{lem:4}$(iv)$ and the construction,
the action of $\sigma $ on~$\Lambda $ swaps $\mathcal A$, $\mathcal B$
and f\/ixes $ I$.
The diagram commutes by these comments
and Lemma~\ref{prop:uhomoh}.
\end{proof}

\begin{Theorem}
\label{thm:faith}
The action of
${\rm  {PSL}}_2(\mathbb Z)$ on $\Delta $ is faithful.
\end{Theorem}

\begin{proof}
Let $g$ denote an element of
${\rm  {PSL}}_2(\mathbb Z)$ that f\/ixes everything in
$\Delta $.
We show that $g=1$.
By Lemma~\ref{lem:conj} and since $\rho$, $\sigma$ generate
${\rm  {PSL}}_2(\mathbb Z)$, there exists an invertible
$\xi \in \Lambda $ such that $g(\eta)=\xi \eta \xi^{-1}$ for
all $\eta  \in \Lambda $. By assumption
$g$ f\/ixes the element $A$ of $\Delta $. Under
the homomomorphim
$\pi : \Delta  \to \Lambda $ the image of $A$ is
$q\mathcal A+q^{-1}{\mathcal A}^{-1}$,
so $g$ f\/ixes
$q\mathcal A+q^{-1}{\mathcal A}^{-1}$ in view of Lemma~\ref{lem:comdiagf}. Therefore $\xi $ commutes
with
$q\mathcal A+q^{-1}{\mathcal A}^{-1}$.
Recall $\mathcal A+{\mathcal A}^{-1} =\mu I$ by
Lemma~\ref{lem:xyz}, so $\xi$ commutes with
$\mathcal A+{\mathcal A}^{-1}$.
By these comments and $q^2\not=1$ we f\/ind
$\xi$ commutes with $\mathcal A$. By a similar
argument $\xi$ commutes with $\mathcal B$. Now
$\xi \in Z(\Lambda )$ by
Lemma~\ref{lem:whatittakes}.
Consequently $g$ f\/ixes everything in $\Lambda $, so $g=1$
 by Lemma~\ref{lem:faith1}.
\end{proof}

\section[A basis for $\Delta $]{A basis for $\boldsymbol{\Delta}$}

 In this section we display a basis for the $\mathbb F$-vector
space $\Delta $.

\begin{Theorem}
\label{thm:basis}
The following is a basis for the $\mathbb F$-vector
space $ \Delta $.
\begin{gather}
\label{eq:basis}
A^iB^jC^k \alpha^r\beta^s\gamma^t, \qquad   i,j,k,r,s,t \geq 0.
\end{gather}
\end{Theorem}

\begin{proof}
We invoke  Bergman's Diamond Lemma
\cite[Theorem~1.2]{berg}.
Consider the symbols
\begin{gather}
\label{eq:letters}
A, \quad B, \quad C, \quad \alpha,\quad \beta, \quad \gamma.
\end{gather}
For an integer $n \geq 0$, by a {\it $\Delta $-word of length $n$}
we mean a sequence
 $x_1x_2\cdots x_n$ such that~$x_i$ is listed in
(\ref{eq:letters}) for $1 \leq i \leq n$.
We interpret the
$\Delta $-word
of length zero to be the multiplicative
identity in $\Delta $.
Consider a
 $\Delta $-word
$w=x_1x_2\cdots x_n$.
By an {\it inversion} for $w$
we mean an ordered pair of integers $(i,j)$
such that
$1\leq i<j\leq n$ and
$x_i$ is strictly to the right of
$x_j$ in the list~(\ref{eq:letters}).
For example $CABA$ has 4 inversions
and
$CB^2A$ has 5 inversions.
A $\Delta $-word is called {\it reducible} whenever it has at least one
inversion, and {\it irreducible} otherwise.
The list
(\ref{eq:basis}) consists of the irreducible $\Delta $-words.
For each integer $n\geq 0$ let
$W_n$ denote the set of $\Delta $-words that have
length~$n$.
Let $W=\cup_{n=0}^\infty W_n$ denote the
set of all
$\Delta $-words.
We now def\/ine a partial
order~$<$ on~$W$.
The def\/inition has two aspects.
$(i)$~For all integers  $n> m\geq 0$,
every word in~$W_m$ is less than
every word in~$W_n$, with respect to $<$.
$(ii)$~For an integer $n\geq 0$
the restriction of $<$ to $W_n$ is
described as follows.
Pick $w, w' \in
 W_n$ and write
 $w=x_1x_2\cdots x_n$.
We say that $w$ {\it covers} $w'$ whenever
there exists an integer~$j$ $(2 \leq j \leq n)$
such that $(j-1,j)$ is an inversion for~$w$,
and
$w'$ is obtained from $w$ by interchanging
$x_{j-1},x_j$.
In this case
$w'$ has one fewer inversions than $w$.
Therefore the transitive closure of the covering
relation on $W_n$ is a partial order on $W_n$, and this is
the restriction of
 $<$ to $W_n$.
We have now def\/ined a partial order $<$ on $W$.
By construction this partial order is
a semi-group partial order
\cite[p.~181]{berg} and satisf\/ies the
descending chain condition
\cite[p.~179]{berg}.
By Def\/inition~\ref{def:uaw} and
Def\/inition~\ref{def:abc}
the def\/ining relations for $\Delta $ can be expressed
as follows:
\begin{gather*}
 BA  =  q^2AB+q\big(q^2-q^{-2}\big)C-q\big(q-q^{-1}\big)\gamma,
\\
 CB  =  q^2BC+q\big(q^2-q^{-2}\big)A-q\big(q-q^{-1}\big)\alpha,
\\
 CA  =  q^{-2}AC+q^{-1}\big(q^{-2}-q^2\big)B-q^{-1}\big(q^{-1}-q\big)\beta,
\\
 \alpha A = A \alpha,
  \qquad
\alpha B = B \alpha,
  \qquad
\alpha C = C \alpha,
\\
 \beta A = A \beta,
  \qquad
\beta B = B \beta,
  \qquad
\beta C = C \beta,
\\
 \gamma A = A \gamma,
  \qquad
\gamma B = B \gamma,
  \qquad
\gamma C = C \gamma,
\\
 \beta \alpha = \alpha \beta,
  \qquad
\gamma \beta = \beta \gamma,
  \qquad
\gamma \alpha = \alpha \gamma.
\end{gather*}
The above equations give reduction rules
for $\Delta$-words,
as we now explain.
Let $w=x_1x_2\cdots x_n$ denote a  reducible $\Delta $-word.
Then there exists an integer $j$ $(2 \leq j \leq n)$
such that $(j-1,j)$ is an inversion for~$w$.
In the above list of equations,
there exists an equation
with $x_{j-1}x_j$ on the left-hand side;
in~$w$ we eliminate $x_{j-1}x_j$ using this equation
and thereby express $w$ as a linear combination of
$\Delta $-words, each less than $w$ with respect to $<$.
Therefore
the reduction rules are compatible
with $<$ in the sense of Bergman \cite[p.~181]{berg}.
In order to employ the Diamond Lemma,
we must show that
the ambiguities
are resolvable in the sense of Bergman
\cite[p.~181]{berg}.
There are potentially two kinds of ambiguities; inclusion
ambiguities and overlap ambiguities
\cite[p.~181]{berg}.
For the present example there are no inclusion ambiguities.
The only nontrivial overlap ambiguity
involves the word
$CBA$. This word can be reduced
in two ways; we could evaluate $CB$ f\/irst or we could
evaluate $BA$ f\/irst.
 Either way, after a three-step reduction
we obtain the same result, which is
\begin{gather*}
 q^{-1}CBA = qABC
+\big(q^2-q^{-2}\big)A^2
-\big(q^2-q^{-2}\big)B^2
+\big(q^2-q^{-2}\big)C^2
\\
\phantom{q^{-1}CBA =}{}
- \big(q-q^{-1}\big)A \alpha
+ \big(q-q^{-1}\big)B \beta
- \big(q-q^{-1}\big)C \gamma .
\end{gather*}
Therefore the
overlap ambiguity $CBA$ is re\-sol\-vable.
We conclude that every  ambiguity is
re\-sol\-vable, so
by the Diamond Lemma
\cite[Theorem~1.2]{berg}
the elements
(\ref{eq:basis}) form a basis for
$\Delta $.
\end{proof}

 On occasion we
wish to discuss the coef\/f\/icients
when an element of $\Delta$ is written as a linear
combination of the elements
(\ref{eq:basis}). To facilitate this discussion
we def\/ine a bilinear form
$(\,,) : \Delta \times \Delta \to \mathbb F$
such that $( u,v) = \delta_{u,v}$
for all elements $u$, $v$ in the basis
(\ref{eq:basis}). In other words
the basis
(\ref{eq:basis}) is
orthonormal with respect to
$( \,,) $.
 Observe that
$(\,,)$ is symmetric.
For $u \in \Delta$,
\begin{gather}
u = \sum \big(u, A^iB^jC^k\alpha^r\beta^s\gamma^t \big)
 A^iB^jC^k\alpha^r\beta^s\gamma^t,
\label{eq:bilsum}
\end{gather}
where the sum is over all elements
 $A^iB^jC^k\alpha^r\beta^s\gamma^t$ in the basis
(\ref{eq:basis}).

\begin{Definition}
\label{def:contrib}
\rm
Let $u \in \Delta $. A given element
 $A^iB^jC^k\alpha^r\beta^s\gamma^t$ in the basis
(\ref{eq:basis}) is said to
{\it contribute to
  $u$} whenever $(u,
 A^iB^jC^k\alpha^r\beta^s\gamma^t )
  \not=0$.
\end{Definition}

\section[A filtration of $\Delta$]{A f\/iltration of $\boldsymbol{\Delta}$}

In this section we obtain
a f\/iltration of $\Delta $
which
is related to the basis from Theorem~\ref{thm:basis}.
 This f\/iltration will be useful
when we investigate the center  $Z(\Delta)$ later in the paper.

We recall some notation.
For subspaces $H$, $K$ of $\Delta $ def\/ine
$HK={\rm Span}\lbrace hk | h \in H, \; k \in K\rbrace$.

\begin{Definition}
\label{def:ui}
We def\/ine subspaces
$\lbrace \Delta_n\rbrace_{n=0}^\infty$
of $\Delta $ such that
\begin{gather*}
\Delta_0=\mathbb F 1,
  \qquad
\Delta_1 =\Delta_0+{\rm Span}\lbrace A,B,C,\alpha,\beta,\gamma\rbrace,
 \qquad
\Delta_n=\Delta_1 \Delta_{n-1}, \qquad n=1,2,\ldots
\end{gather*}
\end{Definition}

\begin{Lemma}
\label{lem:filt}
The following $(i)$--$(iii)$ hold.
\begin{enumerate}\itemsep=0pt
\item[$(i)$]
 $\Delta_{n-1}\subseteq \Delta_n$ for  $n\geq 1$.
\item[$(ii)$]
$\Delta=\cup_{n=0}^\infty \Delta_n$.
\item[$(iii)$]
$\Delta_m \Delta_n = \Delta_{m+n} $ for
 $m,n\geq 0$.
\end{enumerate}
\end{Lemma}

\begin{proof}
$(i)$ Since $\Delta_n=\Delta_1 \Delta_{n-1}$
and $1 \in \Delta_1$.
$(ii)$ Since $A$, $B$, $C$, $\alpha$, $\beta$, $\gamma $ generate~$\Delta$.
$(iii)$ Each side is equal to
$(\Delta_1)^{m+n}$.
\end{proof}

By Lemma~\ref{lem:filt} the sequence
$\lbrace \Delta_n\rbrace_{n=0}^\infty$ is a f\/iltration of~$\Delta$ in the sense of
\cite[p.~202]{carter}.

\begin{Lemma}
\label{lem:almostcom}
Each of the following is contained in $\Delta_1$:
\begin{gather*}
qAB-q^{-1}BA,
\qquad
qBC-q^{-1}CB,
\qquad
qCA-q^{-1}AC.
\end{gather*}
\end{Lemma}

\begin{proof}
Each of the three expressions is a linear combination of
$A$, $B$, $C$, $\alpha$, $\beta$, $\gamma$ and these are contained in~$\Delta_1$.
\end{proof}

\begin{Theorem}\label{lem:unbasis}
For all integers $n\geq 0$ the following is
a basis for the~$\mathbb F$-vector space~$\Delta_n$:
\begin{gather}
\label{eq:basislist}
A^iB^jC^k\alpha^r \beta^s \gamma^t, \qquad
i,j,k,r,s,t\geq 0, \qquad
i+j+k+r+s+t \leq n.
\end{gather}
\end{Theorem}

\begin{proof}
The elements~(\ref{eq:basislist}) are linearly independent by
Theorem
\ref{thm:basis}. We show that the
elements~(\ref{eq:basislist}) span $\Delta_n$.
We will use induction on $n$. Assume $n\geq 2$; otherwise
the result holds by Def\/inition~\ref{def:ui}.
By Def\/inition~\ref{def:ui} $\Delta_n$ is
spanned by the set of elements of the form
$x_1x_2\cdots x_n$
where $x_i \in \lbrace 1, A,B,C,\alpha,\beta,\gamma\rbrace$
for $1 \leq i \leq n$.
Therefore $\Delta_n$ is spanned by the set of elements
of the form
$x_1x_2\cdots x_m$ where $0 \leq m\leq n$ and
$x_i \in \lbrace A,B,C,\alpha,\beta,\gamma\rbrace$
for $1 \leq i \leq m$.
Therefore $\Delta_n$ is spanned by $\Delta_{n-1}$ together with
the set of elements
of the form
$x_1x_2\cdots x_n$ where
$x_i \in \lbrace A,B,C,\alpha,\beta,\gamma\rbrace$
for $1 \leq i \leq n$.
Consider such an element
$x_1x_2\cdots x_n$.
By Lemma~\ref{lem:almostcom} and since
each of $\alpha$, $\beta$, $\gamma$ is central,
we f\/ind that
for $2 \leq j \leq n$,
\begin{gather*}
x_1\cdots  x_{j-1}x_j\cdots  x_n \in
\mathbb F x_1\cdots x_{j}x_{j-1}\cdots  x_n + \Delta_{n-1}.
\end{gather*}
By the above comments
$\Delta_n$ is spanned by $\Delta_{n-1}$ together with the
set
\begin{gather*}
A^iB^jC^k\alpha^r \beta^s \gamma^t, \qquad
i,j,k,r,s,t\geq 0, \qquad
i+j+k+r+s+t = n.
\end{gather*}
By this and induction
$\Delta_n$ is spanned by the elements~(\ref{eq:basislist}).
We have shown that
the elements~(\ref{eq:basislist}) form a basis  for the $\mathbb F$-vector space~$\Delta_n$.
\end{proof}

  Let $V$ denote a vector space over
$\mathbb F$ and let $U$ denote a subspace of~$V$.
By a {\it complement of~$U$ in~$V$} we mean
a subspace $U'$ of $V$ such that
$V=U+U'$ (direct sum).

\begin{Corollary}
\label{lem:com}
For all integers $n\geq 1$ the following is a basis
for a complement of $\Delta_{n-1}$ in $\Delta_n$:
\begin{gather*}
A^iB^jC^k\alpha^r \beta^s \gamma^t, \qquad
i,j,k,r,s,t\geq 0, \qquad
i+j+k+r+s+t = n.
\end{gather*}
\end{Corollary}
\begin{proof}
Use Theorem
\ref{lem:unbasis}.
\end{proof}

\section[The Casimir element $\Omega$]{The Casimir element $\boldsymbol{\Omega}$}

We turn our attention to the center $Z(\Delta)$.
In this section we discuss a certain
 element $\Omega \in Z(\Delta)$
 called the Casimir element.
The name is
motivated by
\cite[equation~(1.3)]{Zhidd}.
 In Section 7
we will use~$\Omega$ to describe
$Z(\Delta)$.
We acknowledge that the results of this section are
extensions of
\cite[Lemma~1]{posta2},
\cite[Section~6]{koro},
\cite[equation~(1.3)]{Zhidd}.

\begin{Lemma}
\label{lem:six}
The following elements of $\Delta $ coincide:
\begin{gather*}
qABC + q^2A^2+q^{-2}B^2+q^2C^2-q A\alpha -q^{-1} B\beta -q C\gamma,
\\
qBCA +
q^2A^2 +
q^2B^2+q^{-2}C^2
-q A\alpha
-qB\beta -q^{-1}C\gamma,
\\
qCAB +
q^{-2}A^2+q^2B^2
+q^2 C^2
-q^{-1} A\alpha -q B \beta
-q C\gamma,
\\
q^{-1}CBA + q^{-2}A^2+q^{2}B^2+q^{-2}C^2-q^{-1}A\alpha -q B\beta -
q^{-1} C\gamma,
\\
q^{-1}ACB +
q^{-2}A^2 +
q^{-2}B^2+q^{2}C^2
-q^{-1} A\alpha
-q^{-1}B\beta -q C\gamma,
\\
q^{-1}BAC +
q^{2}A^2+q^{-2}B^2
+q^{-2}C^2
-q A\alpha
-q^{-1} B \beta
-q^{-1} C\gamma.
\end{gather*}
We denote this common element by $\Omega$.
\end{Lemma}

\begin{proof}
Denote the displayed sequence of elements by
$\Omega^+_B$,
$\Omega^+_C$,
$\Omega^+_A$,
$\Omega^-_B$,
$\Omega^-_C$,
$\Omega^-_A$.
The automorphism
$\rho$ cyclically permutes
$\Omega^+_A$,
$\Omega^+_B$,
$\Omega^+_C$
and cyclically permutes
$\Omega^-_A$, $\Omega^-_B$, $\Omega^-_C$.
The element
$\Omega^+_B-\Omega^-_C$
is
equal to $(q-q^{-1})A$ times
\begin{gather}
\label{eq:need0}
\big(q+q^{-1}\big)A + \frac{qBC-q^{-1}CB}{q-q^{-1}} -\alpha.
\end{gather}
The element
(\ref{eq:need0}) is zero by Def\/inition
\ref{def:abc} so
$\Omega^+_B=\Omega^-_C$.
In this equation we apply $\rho$ twice to get
$\Omega^+_C=\Omega^-_A$
and
$\Omega^+_A=\Omega^-_B$.
The element
$\Omega^+_B-\Omega^-_A$ is equal to
\begin{gather}
\label{eq:need02}
\big(q+q^{-1}\big)C + \frac{qAB-q^{-1}BA}{q-q^{-1}} -\gamma
\end{gather}
times $(q-q^{-1})C$.
The element
(\ref{eq:need02}) is zero by Def\/inition
\ref{def:abc} so
$\Omega^+_B=\Omega^-_A$.
Applying $\rho$ twice we get
$\Omega^+_C=\Omega^-_B$
and
$\Omega^+_A=\Omega^-_C$.
By these comments
$\Omega^+_B$,
$\Omega^+_C$,
$\Omega^+_A$,
$\Omega^-_B$,
$\Omega^-_C$,
$\Omega^-_A$ coincide.
\end{proof}

\begin{Theorem}
\label{lem:center}
The element $\Omega$ from
Lemma~{\rm \ref{lem:six}} is central in $\Delta $.
\end{Theorem}

\begin{proof}
We f\/irst show $\Omega A = A \Omega $.
We will work with the equations~(\ref{eq:u2}),
(\ref{eq:u3}) from Def\/inition~\ref{def:abc}. Consider the equation which is
$qC $ times~(\ref{eq:u2}) plus~(\ref{eq:u2}) times $q^{-1}C$
minus $\gamma$ times~(\ref{eq:u2})
plus $\beta$ times~(\ref{eq:u3})
minus
$q^{-1}B$ times~(\ref{eq:u3})
minus~(\ref{eq:u3}) times~$qB$.
After some cancellation this equation yields
$\Omega^+_B A- A\Omega^+_C =0$, where
$\Omega^+_B$, $\Omega^+_C$
are from the proof of Lemma
\ref{lem:six}.
Therefore
$\Omega A= A\Omega $.
 One similarly
f\/inds
$
\Omega B = B \Omega$ and
$
\Omega C = C \Omega$.
The elements $A$, $B$, $C$ generate~$\Delta $ so~$\Omega$ is central in~$\Delta$.
\end{proof}

\begin{Definition}
We call $\Omega$ the {\it Casimir} element
of~$\Delta $.
\end{Definition}

\begin{Theorem}
The Casimir element $\Omega$ is fixed by everything
in ${\rm  {PSL}}_2(\mathbb Z)$.
\end{Theorem}

\begin{proof}
Since $\rho$, $\sigma$ generate
${\rm  {PSL}}_2(\mathbb Z)$ it suf\/f\/ices to show
that each of $\rho$, $\sigma$ f\/ixes $\Omega$.
We use the notation $\Omega^+_A$, $\Omega^+_B$
from the proof of Lemma~\ref{lem:six}.
Observe that $\rho$ f\/ixes $\Omega$ since
$\rho(\Omega^+_A)=\Omega^+_B$.
To verify that  $\sigma$ f\/ixes $\Omega$ we show
that $\sigma(\Omega^+_B)=\Omega^+_A$.
For notational convenience def\/ine
\begin{gather*}
C' = C + \frac{AB-BA}{q-q^{-1}}.
\end{gather*}
By Theorem~\ref{thm:try} and the def\/inition~$\Omega^+_B$,
\begin{gather*}
\sigma(\Omega^+_B) =
qBAC' + q^2B^2+q^{-2}A^2+q^2C'^2-q B\beta-q^{-1}A \alpha -q C' \gamma .
\end{gather*}
By this and the def\/inition of $\Omega^+_A$,
\begin{gather*}
\sigma(\Omega^+_B)-\Omega^+_A
=
\big(BA+q^{-1}C+qC'-\gamma\big)qC'
 - qC\big(AB+qC+q^{-1}C'-\gamma\big).
\end{gather*}
In the above equation each parenthetical expression
is zero
so
$\sigma(\Omega^+_B)=\Omega^+_A$.
Therefore~$\sigma$ f\/ixes~$\Omega$.
\end{proof}

\section[A basis for~$\Delta $ that involves~$\Omega$]{A basis for $\boldsymbol{\Delta}$ that involves~$\boldsymbol{\Omega}$}

 In Theorem~\ref{thm:basis} we displayed a basis for
$\Delta $. In this section we display a related
basis for $\Delta $ that involves the Casimir element~$\Omega$.
In the next section we will use the related basis
to describe  the center~$Z(\Delta )$.

 Recall the f\/iltration $\lbrace \Delta_n \rbrace_{n=0}^\infty$
of $\Delta $ from Def\/inition~\ref{def:ui}.

\begin{Lemma}
\label{lem:omegapower}
For all integers $\ell \geq 1$
the following hold:
\begin{enumerate}\itemsep=0pt
\item[$(i)$] $\Omega^\ell \in \Delta_{3 \ell}$.
\item[$(ii)$] $\Omega^\ell -q^{\ell^2}A^\ell B^\ell C^\ell  \in
\Delta_{3\ell-1}$.
\end{enumerate}
\end{Lemma}

\begin{proof}
Consider the expression for $\Omega$ from the
f\/irst displayed line of
Lemma~\ref{lem:six}.
The term $qABC$ is in $\Delta_3$
and the remaining terms are in $\Delta_2$. Therefore
$\Omega \in \Delta_3$ and $\Omega-qABC \in \Delta_2$. By this and
Lemma~\ref{lem:filt}$(iii)$ we f\/ind
$\Omega^\ell \in \Delta_{3\ell}$ and
$\Omega^\ell - q^\ell (ABC)^\ell
\in \Delta_{3\ell-1}$.
Using Lemma~\ref{lem:almostcom}  we obtain
$(ABC)^\ell -q^{\ell (\ell-1)}A^\ell B^\ell C^\ell \in \Delta_{3\ell -1}$.
By these comments
$\Omega^\ell -q^{\ell^2}A^\ell B^\ell C^\ell  \in \Delta_{3\ell-1}$.
\end{proof}

\begin{Lemma}
\label{lem:basis2c}
For all integers $n\geq 1$
the following is a basis
for a complement of
$\Delta_{n-1}$ in $\Delta_n$:
\begin{gather*}
 A^iB^jC^k\Omega^\ell \alpha^r \beta^s \gamma^t, \qquad
i,j,k,\ell, r,s,t\geq 0,
\qquad ijk=0,
 \qquad
i+j+k+3\ell+r+s+t = n.
\end{gather*}
\end{Lemma}

\begin{proof}
Let $\mathbb I_n$ denote the set consisting of the 6-tuples
$(i,j,k,r,s,t)$ of nonnegative integers whose sum is~$n$.
By Corollary~\ref{lem:com} the
following
is a basis for a complement of
$\Delta_{n-1}$ in $\Delta_n$:
\begin{gather*}
A^iB^jC^k\alpha^r\beta^s\gamma^t, \qquad
(i,j,k,r,s,t) \in {\mathbb I}_n.
\end{gather*}
Let $\mathbb J_n$ denote the set consisting of the 7-tuples
$(i,j,k,\ell,r,s,t)$ of nonnegative integers such that
$ijk=0$ and
$i+j+k+3\ell+r+s+t = n$.
Observe that the map
\begin{gather*}
\mathbb J_n  \to  \mathbb I_n,
\qquad
(i,j,k,\ell,r,s,t)  \mapsto   (i+\ell,j+\ell,k+\ell,r,s,t)
\end{gather*}
is a bijection.
Suppose we are given $(i,j,k,\ell,r,s,t) \in \mathbb J_n$.
By
Lemma \ref{lem:omegapower},
$\Delta_{n-1}$ contains
\begin{gather*}
A^iB^jC^k\Omega^\ell \alpha^r \beta^s \gamma^t -
q^{\ell^2} A^iB^jC^k A^\ell B^\ell C^\ell \alpha^r \beta^s \gamma^t.
\end{gather*}
By Lemma~\ref{lem:almostcom},  $\Delta_{n-1}$ contains
\begin{gather*}
A^iB^jC^k A^\ell B^\ell C^\ell \alpha^r \beta^s \gamma^t
-
q^{2j \ell} A^{i+\ell}B^{j+\ell}C^{k+\ell}  \alpha^r \beta^s \gamma^t.
\end{gather*}
Therefore $\Delta_{n-1}$ contains
\begin{gather*}
A^iB^jC^k\Omega^\ell \alpha^r \beta^s \gamma^t -
q^{\ell(2j+ \ell)} A^{i+\ell}B^{j+\ell}C^{k+\ell}  \alpha^r \beta^s \gamma^t.
\end{gather*}
By these comments the following is a basis for a complement of
$\Delta_{n-1}$ in $\Delta_n$:
\begin{gather*}
A^iB^jC^k\Omega^\ell \alpha^r\beta^s\gamma^t, \qquad
(i,j,k,\ell,r,s,t) \in \mathbb J_n.
\end{gather*}
The result follows.
\end{proof}

\begin{Note}
Pick an integer $n \geq 1$.
In
Corollary
\ref{lem:com} and
Lemma
\ref{lem:basis2c}
we mentioned a
complement of $\Delta_{n-1}$ in $\Delta_n$. These
complements are not the same
in general.
\end{Note}

\begin{Proposition}
\label{lem:basis2n}
For all integers $n\geq 0$
the following is a basis
for
the $\mathbb F$-vector space
$\Delta_n$:
\begin{gather*}
 A^iB^jC^k\Omega^\ell \alpha^r \beta^s \gamma^t, \qquad
i,j,k,\ell, r,s,t\geq 0,
\qquad ijk=0,
 \qquad
i+j+k+3\ell+r+s+t \leq  n.
\end{gather*}
\end{Proposition}

\begin{proof}
By Lemma~\ref{lem:basis2c} and
$\Delta_0=\mathbb F 1$.
\end{proof}

\begin{Theorem}
\label{lem:basis2}
The following is a basis for the $\mathbb F$-vector space
$\Delta$:
\begin{gather*}
A^iB^jC^k\Omega^\ell \alpha^r \beta^s \gamma^t, \qquad
i,j,k,\ell, r,s,t\geq 0,
 \qquad ijk=0.
\end{gather*}
\end{Theorem}

\begin{proof}
Combine
Lemma~\ref{lem:filt}$(ii)$
and
Proposition~\ref{lem:basis2n}.
\end{proof}

\section[The center $Z(\Delta)$]{The center $\boldsymbol{Z(\Delta)}$}

In this section we give a detailed description of
the center $Z(\Delta)$, under the assumption that $q$ is not a root of
 unity.
For such $q$ we show that $Z(\Delta)$ is generated by
$\Omega$, $\alpha$, $\beta$, $\gamma$ and
 isomorphic to a polynomial algebra in four
variables.

  Recall the commutator
 $\lbrack r,s\rbrack=rs-sr$.

\begin{Lemma}
\label{lem:ijk}
Let
$i$, $j$, $k$
denote nonnegative integers.
Then $\Delta_{i+j+k}$ contains each of the fol\-lo\-wing:
\begin{gather}
\label{eq:1}
 \big\lbrack A,A^iB^jC^k\big\rbrack - \big(1-q^{2j-2k}\big)A^{i+1}B^jC^k,
\\
\label{eq:2}
 \big\lbrack B,A^iB^jC^k\big\rbrack - \big(q^{2i}-q^{2k}\big)A^iB^{j+1}C^k,
\\
\label{eq:3}
 \big\lbrack C,A^iB^jC^k\big\rbrack - \big(q^{2j-2i}-1\big)A^iB^jC^{k+1}.
\end{gather}
\end{Lemma}

\begin{proof}
Concerning (\ref{eq:1}), observe
\begin{gather*}
\big\lbrack A,A^iB^jC^k\big\rbrack = A^{i+1}B^jC^k -A^iB^jC^kA.
\end{gather*}
By Lemma
\ref{lem:almostcom} $\Delta_{i+j+k}$ contains
\begin{gather*}
A^iB^jC^kA - q^{2j-2k}A^{i+1}B^jC^k.
\end{gather*}
By these comments
$\Delta_{i+j+k}$ contains
(\ref{eq:1}).
One similarly f\/inds that
$\Delta_{i+j+k}$ contains
(\ref{eq:2}),
(\ref{eq:3}).
\end{proof}

\begin{Theorem}
\label{thm:centerbasis}
Assume that $q$ is not a root of unity. Then the following is a basis for
the $\mathbb F$-vector space~$Z(\Delta)$.
\begin{gather}
\Omega^\ell \alpha^r \beta^s\gamma^t, \qquad  \ell,r,s,t \geq 0.
\label{eq:cbasis}
\end{gather}
\end{Theorem}

\begin{proof}
Abbreviate $Z=Z(\Delta)$.
The elements~(\ref{eq:cbasis})
are linearly independent by
Theorem~\ref{lem:basis2},
so it suf\/f\/ices to
show that they span~$Z$.
Let~$Z'$ denote the subspace of~$\Delta$ spanned by~(\ref{eq:cbasis}), and note that
$Z' \subseteq Z$. To show
$Z'= Z$, we assume that $Z'$ is properly contained in
$Z$ and get a contradiction.
Def\/ine the set $E=Z\backslash Z'$
and note $E\not=\varnothing$.
We have
$\Delta_0=\mathbb F 1
\subseteq Z'$ so
$E \cap \Delta_0 =\varnothing$.
By this and
Lemma~\ref{lem:filt}$(i)$,$(ii)$ there exists a unique integer
$n\geq 1$ such that
$E\cap \Delta_n \not=\varnothing$ and
$E\cap \Delta_{n-1} =\varnothing$.
Fix $u \in E\cap \Delta_n$.
Let
$S=S(u)$ denote
the set of 6-tuples $(i,j,k,r,s,t)$ of nonnegative integers
whose sum is $n$ and
 $A^iB^jC^k\alpha^r\beta^s\gamma^t$
contributes to $u$ in
the sense of
Def\/inition~\ref{def:contrib}.
By~(\ref{eq:bilsum}) and Corollary~\ref{lem:com},
$\Delta_{n-1}$ contains
\begin{gather}
\label{eq:dif}
u- \sum_{(i,j,k,r,s,t) \in S}
\big( u,
A^iB^jC^k\alpha^r \beta^s \gamma^t \big)
A^iB^jC^k\alpha^r \beta^s \gamma^t.
\end{gather}
We are going to show that
$i=j=k$ for all
$(i,j,k,r,s,t) \in S$.
To this end
we f\/irst claim that $j=k$ for all
$(i,j,k,r,s,t) \in S$.
To prove the claim,
take the commutator of $A$ with~(\ref{eq:dif}) and evaluate the
result using the following facts.
By construction $u \in E \subseteq Z$
so
  $\lbrack A,u\rbrack=0$.
By Lemma~\ref{lem:filt}$(iii)$
$A\Delta_{n-1} \subseteq  \Delta_n $
and
$\Delta_{n-1}A \subseteq \Delta_n$
so
 $\lbrack A,\Delta_{n-1}\rbrack \subseteq \Delta_n$.
Moreover
 each of~$\alpha$,~$\beta$,~$\gamma$
 is central.
The above evaluation shows that~$\Delta_n$ contains
\begin{gather*}
\sum_{(i,j,k,r,s,t) \in S}
\big( u, A^iB^jC^k
\alpha^r \beta^s \gamma^t \big)
\big\lbrack A, A^iB^jC^k \big\rbrack
\alpha^r \beta^s \gamma^t.
\end{gather*}
Pick
$(i,j,k,r,s,t) \in S$.
By
Lemma~\ref{lem:filt}$(iii)$ and
Lemma~\ref{lem:ijk}, $\Delta_n$ contains
\begin{gather*}
\big\lbrack A, A^iB^jC^k \big\rbrack
\alpha^r \beta^s \gamma^t
-
\big(1-q^{2j-2k}\big)  A^{i+1}B^jC^k
\alpha^r \beta^s \gamma^t.
\end{gather*}
By the above comments
$\Delta_n$ contains
\begin{gather*}
\sum_{
(i,j,k,r,s,t) \in S}
\big(u, A^iB^jC^k
\alpha^r \beta^s \gamma^t \big)
\big(1-q^{2j-2k}\big)  A^{i+1}B^jC^k
\alpha^r \beta^s \gamma^t.
\end{gather*}
For all
$(i,j,k,r,s,t) \in S$
the element
$A^{i+1}B^jC^k
\alpha^r \beta^s \gamma^t$
is contained
in the basis for the complement of $\Delta_n$ in $\Delta_{n+1}$
given in Corollary~\ref{lem:com}. Therefore
\begin{gather*}
\big(u, A^iB^jC^k
\alpha^r \beta^s \gamma^t\big)
\big(1-q^{2j-2k}\big)  = 0
\qquad
\forall \;(i,j,k,r,s,t) \in S.
\end{gather*}
By the def\/inition of $S$ we have
$( u, A^iB^jC^k
\alpha^r \beta^s \gamma^t )\not=0$
 for
all
$(i,j,k,r,s,t) \in S$.
Therefore
$1-q^{2j-2k} = 0 $ for all
$(i,j,k,r,s,t) \in S$.
The scalar $q$ is not a root of unity
so $j=k$ for all
$(i,j,k,r,s,t) \in S$.
The claim is proved.
We next claim that
$i=j$ for all
$(i,j,k,r,s,t) \in S$.
This claim is proved like the previous one,
 except that
as we begin the argument below~(\ref{eq:dif}),
we use~$C$ instead of~$A$ in the commutator.
By the two claims
$i=j=k$ for all
$(i,j,k,r,s,t) \in S$.
In light of this we revisit the assertion above~(\ref{eq:dif}) and conclude that
$\Delta_{n-1}$ contains
\begin{gather*}
u-\sum_{
(i,i,i,r,s,t) \in S}
\big(u, A^iB^iC^i
\alpha^r \beta^s \gamma^t \big)
A^{i}B^iC^i
\alpha^r \beta^s \gamma^t.
\end{gather*}
Pick
$(i,i,i,r,s,t) \in S$.
By
Lemma~\ref{lem:filt}$(iii)$ and
Lemma~\ref{lem:omegapower},
$\Delta_{n-1}$ contains
\begin{gather*}
A^iB^iC^i
\alpha^r \beta^s \gamma^t
-q^{-i^2}\Omega^i \alpha^r \beta^s \gamma^t.
\end{gather*}
By these comments
$\Delta_{n-1}$ contains
\begin{gather*}
\label{eq:summ}
u-
\sum_{
(i,i,i,r,s,t) \in S}
q^{-i^2}
\big(u, A^iB^iC^i
\alpha^r \beta^s \gamma^t  \big)
\Omega^i \alpha^r \beta^s \gamma^t.
\end{gather*}
In the above expression let
$\psi$ denote the main sum,
so that
$u-\psi \in \Delta_{n-1}$.
Observe
$\psi \in Z' \subseteq Z$.
Recall $u \in E=Z\backslash Z'$ so
$u \in Z$ and $u \not\in Z'$. By these comments
$u - \psi \in Z$ and
$u - \psi \not\in Z'$.
Therefore $u - \psi \in E$ so
 $u - \psi \in E\cap \Delta_{n-1}$.
This contradicts
 $E\cap \Delta_{n-1}=\varnothing $ so
$Z=Z'$. The result follows.
\end{proof}

 We mention two corollaries of Theorem~\ref{thm:centerbasis}.

\begin{Corollary}
Assume that $q$ is not a root of unity. Then
$Z(\Delta)$ is generated by $\Omega$, $\alpha$, $\beta$, $\gamma$.
\end{Corollary}

  Let $\lbrace \lambda_i\rbrace_{i=1}^4$ denote mutually
commuting indeterminates. Let
$\mathbb F\lbrack \lambda_1,\lambda_2,\lambda_3,\lambda_4\rbrack$
denote the $\mathbb F$-algebra consisting of the polynomials
in
$\lbrace \lambda_i\rbrace_{i=1}^4$
that have all coef\/f\/icients in~$\mathbb F$.

\begin{Corollary}
Assume that $q$ is not a root of unity. Then
there exists an $\mathbb F$-algebra isomorphism
$Z(\Delta)\to
\mathbb F\lbrack \lambda_1,\lambda_2,\lambda_3,\lambda_4\rbrack$
that sends
\begin{gather*}
\Omega \mapsto \lambda_1,
\qquad
\alpha \mapsto \lambda_2,
\qquad
\beta \mapsto \lambda_3,
\qquad
\gamma \mapsto \lambda_4.
\end{gather*}
\end{Corollary}

\section[The $q$-Onsager algebra $\mathcal O$]{The $\boldsymbol{q}$-Onsager algebra $\boldsymbol{\mathcal O}$}

 In the theory of tridiagonal pairs
there is an algebra known as the tridiagonal algebra
\cite[Def\/inition~3.9]{qSerre}. This algebra
is def\/ined using several parameters, and for
a certain value of these parameters
the algebra
is sometimes called the $q$-Onsager algebra
 $\mathcal O$
\cite[Section~4]{basnc}.
Our next goal is to show how
 $\mathcal O$ and
 $\Delta$ are related.
In this section we def\/ine
 $\mathcal O$ and discuss some of its properties.
In the next section
we will relate
 $\mathcal O$ and $\Delta$.

\begin{Definition}[\protect{\cite[Def\/inition~3.9]{qSerre}}]
\label{def:qons}
Let ${\mathcal O}={\mathcal O}_q$ denote the
 $\mathbb F$-algebra
 def\/ined by genera\-tors~$X$,~$Y$ and relations
 \begin{gather}
X^3Y-\lbrack 3\rbrack_q X^2YX+\lbrack 3\rbrack_q XYX^2-YX^3
= -\big(q^2-q^{-2}\big)^2(XY-YX),
\label{eq:td1}
\\
Y^3X-\lbrack 3\rbrack_q Y^2XY+\lbrack 3\rbrack_q YXY^2-XY^3
= -\big(q^2-q^{-2}\big)^2(YX-XY).
\label{eq:td2}
 \end{gather}
We call $\mathcal O$ the {\it $q$-Onsager algebra}.
\end{Definition}

  The following def\/inition
is motivated by
Theorem~\ref{prop:altpres}.

\begin{Definition}
\label{lem:RS}
Let $\xi_1$, $\xi_2$ denote the following elements of
$\mathcal O$:
\begin{gather*}
\xi_1 =  XY-YX,
\\
 \xi_2 =
X^2Y^2-Y^2X^2+\big(q^2+q^{-2}\big)(YXYX-XYXY).
\end{gather*}
\end{Definition}

Referring to Def\/inition~\ref{lem:RS},
we are going to show that
$\xi_1$, $\xi_2$
commute if and only if~$q^6 \not=1$.
We will use the following results, which
apply to any $\mathbb F$-algebra.

\begin{Lemma}
\label{lem:gen}
Let $x$, $y$ denote elements in any $\mathbb F$-algebra,
and consider the commutator
\begin{gather}
\big\lbrack
xy-yx, x^2y^2-y^2x^2 +\big(q^2+q^{-2}\big)(yxyx-xyxy)
\big\rbrack.
\label{lem:com12}
\end{gather}
\begin{enumerate}\itemsep=0pt
\item[$(i)$] The element
\eqref{lem:com12} is equal to
\begin{gather*}
 xyx^2y^2-x^2y^2xy+yxy^2x^2-y^2x^2yx
+
x^2y^3x-xy^3x^2+y^2x^3y-yx^3y^2
\\
 \qquad {} -
\big(q^2+q^{-2}\big)\big(
xyxy^2x
-xy^2xyx
+yxyx^2y
-yx^2yxy
\big).
\end{gather*}
\item[$(ii)$] The element
\eqref{lem:com12} times
$\lbrack 3\rbrack_q$ is equal to
\begin{gather}
\lbrack y,
\lbrack y,
\lbrack x,
\lbrack x,
\lbrack x,
y
\rbrack \,
\rbrack_q
\rbrack_{q^{-1}}
\rbrack_q
\rbrack_{q^{-1}}
+
\lbrack x,
\lbrack x,
\lbrack y,
\lbrack y,
\lbrack y,
x
\rbrack \,
\rbrack_q
\rbrack_{q^{-1}}
\rbrack_q
\rbrack_{q^{-1}},
\label{eq:double}
\end{gather}
where $\lbrack u,v\rbrack_\epsilon $ means $\epsilon uv-
\epsilon^{-1}vu$.
\end{enumerate}
\end{Lemma}

\begin{proof}
$(i)$ Expand
(\ref{lem:com12}) and simplify the result.
$(ii)$
Expand
(\ref{eq:double})
and compare it with the expression in $(i)$ above.
\end{proof}

We return our attention to the
elements $\xi_1$, $\xi_2$ in $\mathcal O$.

\begin{Proposition}
\label{lem:RS2}
Referring to
Definition {\rm \ref{lem:RS}},
the elements $\xi_1$, $\xi_2$ commute
if and only if
$q^6 \not=1$.
\end{Proposition}

\begin{proof}
First assume
$q^6 \not=1$, so that
$\lbrack 3 \rbrack_q$ is nonzero.
Applying
Lemma
\ref{lem:gen}$(ii)$ to the elements $x=X$ and $y=Y$ in
the algebra $\mathcal O$, we f\/ind that
$\lbrack \xi_1, \xi_2\rbrack $ times $\lbrack 3\rbrack_q$
is equal to
\begin{gather}
\lbrack Y,
\lbrack Y,
\lbrack X,
\lbrack X,
\lbrack X,
Y
\rbrack \,
\rbrack_q
\rbrack_{q^{-1}}
\rbrack_q
\rbrack_{q^{-1}}
+
\lbrack X,
\lbrack X,
\lbrack Y,
\lbrack Y,
\lbrack Y,
X
\rbrack \,
\rbrack_q
\rbrack_{q^{-1}}
\rbrack_q
\rbrack_{q^{-1}}.
\label{eq:bigf}
\end{gather}
We show that the element (\ref{eq:bigf}) is zero.
Observe
\begin{gather*}
\lbrack X,
\lbrack X,
\lbrack X,
Y
\rbrack \,
\rbrack_q
\rbrack_{q^{-1}}
 =
X^3Y-\lbrack 3\rbrack_q X^2YX+\lbrack 3\rbrack_q XYX^2-YX^3
\\
\phantom{\lbrack X,
\lbrack X,
\lbrack X,
Y
\rbrack \,
\rbrack_q
\rbrack_{q^{-1}}}{} = -\big(q^2-q^{-2}\big)^2 \lbrack X,Y \rbrack
= \big(q^2-q^{-2}\big)^2 \lbrack Y,X \rbrack.
\end{gather*}
Similarly
\begin{gather*}
\lbrack Y,
\lbrack Y,
\lbrack Y,
X
\rbrack \,
\rbrack_q
\rbrack_{q^{-1}}
= \big(q^2-q^{-2}\big)^2 \lbrack X,Y \rbrack.
\end{gather*}
Therefore
\begin{gather*}
\lbrack Y,
\lbrack Y,
\lbrack X,
\lbrack X,
\lbrack X,
Y
\rbrack \,
\rbrack_q
\rbrack_{q^{-1}}
\rbrack_q
\rbrack_{q^{-1}}
 =
\big(q^2-q^{-2}\big)^2\lbrack Y,
\lbrack Y,
\lbrack Y,
X
\rbrack \,
\rbrack_q
\rbrack_{q^{-1}}
=
(q^2-q^{-2})^4\lbrack X,Y\rbrack.
\end{gather*}
Similarly
\begin{gather*}
\lbrack X,
\lbrack X,
\lbrack Y,
\lbrack Y,
\lbrack Y,
X
\rbrack \,
\rbrack_q
\rbrack_{q^{-1}}
\rbrack_q
\rbrack_{q^{-1}}
 =
\big(q^2-q^{-2}\big)^4\lbrack Y,X\rbrack .
\end{gather*}
By these comments
the element (\ref{eq:bigf})
is zero, and therefore
$\xi_1$, $\xi_2$ commute.

  Next assume
$q^6=1$, so that
$\lbrack 3\rbrack_q=0$ and
$(q^2-q^{-2})^2=-3$.
In this case the relations
(\ref{eq:td1}),
(\ref{eq:td2}) become
\begin{gather*}
X^3Y-YX^3=3(XY-YX),
\qquad
Y^3X-XY^3=3(YX-XY).
\end{gather*}
In the above line the relation on the left (resp. right)
asserts that $X^3-3X$ (resp.~$Y^3-3Y$)
commutes with $Y$ (resp.~$X$) and is therefore
central in~$\mathcal O$.
We show that $\xi_1$, $\xi_2$ do not commute by
displaying
an $\mathcal O$-module on which
$\lbrack \xi_1, \xi_2\rbrack$ is nonzero.
Def\/ine a group $G$ to be the free product
$H * K$ where each of~$H$,~$K$ is a cylic
group of order~3.
Let $h$ (resp.~$k$) denote a generator for~$H$ (resp.~$K$).
Let $\mathbb F G$ denote the group $\mathbb F$-algebra.
We give $\mathbb F G$ an $\mathcal O$-module structure.
To do this we specify the action of~$X$,~$Y$ on
the basis $G$ of $\mathbb F G$. For the identity $ 1\in G$
let $X.1=h$ and $Y.1=k$.
For $1 \not=g \in G$ write
$g=g_1g_2\cdots g_n$ such that
 for $1 \leq i \leq n$,
 $ g_i \in \lbrace h, h^{-1}\rbrace $
 for one parity of $i$ and
$g_i \in \lbrace k,k^{-1} \rbrace$ for the other parity of~$i$.
Let~$X$,~$Y$ act on~$g$ as follows:
\begin{center}
\begin{tabular}{c|c c c c}
Case &  $g_1=h$ & $g_1=h^{-1}$  & $g_1=k$ & $g_1=k^{-1}$ \\
\hline \hline
$X.g$ & $hg$ & $ 3h^{-1}g$ & $hg$ & $hg$ \\
$Y.g$ & $kg$ & $ kg$ & $kg$ & $3k^{-1}g$
    \end{tabular}
        \end{center}
We have now specif\/ied the actions of $X$, $Y$
on~$G$.
These actions
give $\mathbb F G$ an
 $\mathcal O$-module
structure on
which $X^3=3X$ and $Y^3=3Y$.
For the $\mathcal O$-module $\mathbb F G$ we now
apply
$\lbrack \xi_1, \xi_2\rbrack$ to the vector $1$.
By
Lemma~\ref{lem:gen}$(i)$ and the construction,
the element
$\lbrack \xi_1, \xi_2\rbrack$ sends  $1$
to
\begin{gather*}
 hkh^{-1}k^{-1}-h^{-1}k^{-1}hk+
khk^{-1}h^{-1}-k^{-1}h^{-1}kh
+
  3 h^{-1}kh-3hkh^{-1}
+ 3 k^{-1}hk-3khk^{-1}
\\
\qquad {} +
hkhk^{-1}h-hk^{-1}hkh+
khkh^{-1}k-kh^{-1}khk.
\end{gather*}
The above element of $\mathbb F G$
is nonzero. Therefore
$\lbrack \xi_1, \xi_2\rbrack \not=0$
so
$\xi_1$, $\xi_2$ do not commute.
\end{proof}

\section[How $\mathcal O$ and $\Delta$ are related]{How $\boldsymbol{\mathcal O}$ and $\boldsymbol{\Delta}$ are related}

 Recall the $q$-Onsager algebra $\mathcal O$
from Def\/inition~\ref{def:qons}. In this section we discuss how~$\mathcal O$ and~$\Delta$ are related.

\begin{Definition}  Let $\lambda $ denote
an indeterminate that commutes with
everything in $\mathcal O$.
Let~${\mathcal O}\lbrack \lambda \rbrack$
denote the $\mathbb F$-algebra consisting of
the polynomials in $\lambda $ that have all coef\/f\/icients
in~$\mathcal O$.
We view~$\mathcal O$ as an $\mathbb F$-subalgebra of
 ${\mathcal O}\lbrack \lambda \rbrack$.
\end{Definition}

\begin{Lemma}
\label{lem:varphi}
There exists a unique
$\mathbb F$-algebra
homomorphism $ \varphi: {\mathcal O}\lbrack \lambda \rbrack \to \Delta $
that sends
\begin{gather*}
X\to A, \qquad
Y\to B, \qquad
\lambda \to \gamma.
\end{gather*}
Moreover $\varphi$ is surjective.
\end{Lemma}

\begin{proof}
Compare
Theorem~\ref{prop:altpres} and
Def\/inition~\ref{def:qons}.
\end{proof}

Let $J$ denote the 2-sided ideal of
 ${\mathcal O}\lbrack \lambda \rbrack $
 generated by
$(q-q^{-1})^2\xi_1 \lambda +\xi_2$,
where $\xi_1$, $\xi_2$ are from
Def\/inition~\ref{lem:RS}.
By Theorem~\ref{prop:altpres},
$J$ is the
kernel of~$ \varphi$.
Therefore~$\varphi$ induces an isomorphism of
$\mathbb F$-algebras
${\mathcal O}\lbrack \lambda \rbrack /J
\to \Delta $.

\begin{Theorem}
\label{thm:97}
The $\mathbb F$-algebra $\Delta $ is isomorphic to
${\mathcal O}\lbrack \lambda \rbrack /J$,
where~$J$ is the $2$-sided ideal of
 ${\mathcal O}\lbrack \lambda \rbrack$
 generated by
$(q-q^{-1})^2\xi_1 \lambda +\xi_2$.
\end{Theorem}

 We now adjust our point of view.

\begin{Definition}
\label{def:hatvp}
Let $\phi: {\mathcal O} \to \Delta $
denote the $\mathbb F$-algebra homomorphism that sends
$X\mapsto A$ and $Y\mapsto B$. Observe that
$\phi$ is the restriction of
$ \varphi $ to~${\mathcal O}$.
\end{Definition}

 We now describe the image and kernel of the
 homomorphism $\phi$ in Def\/inition~\ref{def:hatvp}.
We will use the following notation.

\begin{Definition}
\label{def:mathcalS}
For any subset $\mathcal S \subseteq \Delta$
let $\langle \mathcal S\rangle$ denote the
$\mathbb F$-subalgebra of $\Delta$ generated by
$\mathcal S$.
\end{Definition}

\begin{Lemma}
\label{def:hatvpbasic}
For the homomorphism
$\phi :\mathcal O \to \Delta $ from
Definition~{\rm \ref{def:hatvp}},
the image is
$\langle A,B\rangle$. The kernel is
${\mathcal O} \cap J$,
where $J$ is defined above
Theorem~{\rm \ref{thm:97}}.
\end{Lemma}

\begin{proof}
Routine.
\end{proof}

We are going to show
that
$\phi$ is not injective. To do this
we display some nonzero elements in the
kernel~${\mathcal O} \cap J$.

\begin{Lemma}
\label{lem:RS3}
${\mathcal O} \cap J$
contains
the elements
\begin{gather*}
\xi_1 z  \xi_2-\xi_2 z \xi_1,
\qquad
z \in {\mathcal O}.
\end{gather*}
Here $\xi_1$, $\xi_2$ are
from Definition~{\rm \ref{lem:RS}}.
\end{Lemma}

\begin{proof}
Let $z $ be given.
Each of
$z$, $\xi_1$, $\xi_2$ is in
$\mathcal O$ so
$\xi_1 z  \xi_2-\xi_2 z \xi_1 \in \mathcal O$.
The element
$\xi_1 z  \xi_2-\xi_2 z \xi_1$ is equal to
\begin{gather*}
\xi_1 z\big( \big(q-q^{-1}\big)^2 \xi_1 \lambda +  \xi_2\big)
-
\big(\big(q-q^{-1}\big)^2 \xi_1 \lambda +  \xi_2\big)
z \xi_1
\end{gather*}
and is therefore contained in $J$. The result follows.
\end{proof}

We now display
some elements $z$ in $\mathcal O$ such that
$\xi_1 z \xi_2-\xi_2 z \xi_1$ is nonzero.
For general $q$ we cannot take $z=1$ in view of
Proposition~\ref{lem:RS2}, so we proceed to the next simplest
case.

\begin{Lemma}
\label{lem:nonz}
The following  elements of~$\mathcal O$ are nonzero:
\begin{gather*}
\xi_1 X \xi_2-\xi_2 X \xi_1,
\qquad
\xi_1 Y \xi_2-\xi_2 Y \xi_1.
\end{gather*}
Moreover
${\mathcal O} \cap J$ is nonzero.
\end{Lemma}

\begin{proof}
To show $\xi_1 X \xi_2-\xi_2 X \xi_1$ is nonzero
we display an $\mathcal O$-module
on which
$\xi_1 X \xi_2-\xi_2 X \xi_1$ is nonzero.
This $\mathcal O$-module is a variation on an $\mathcal O$-module
due to M.~Vidar
\cite[Theorem 9.1]{Vidar}.
Def\/ine $\theta_i = q^{2i}+q^{-2i}$
for $i=0,1,2$. Def\/ine $\vartheta = (q^4-q^{-4})(q^2-q^{-2})(q-q^{-1})^2$.
Adapting~\cite[Theorem 9.1]{Vidar}
there exists a four-dimensional $\mathcal O$-module $V$
with the following property: $V$ has a~basis with respect to
which the matrices representing~$X$, $Y$ are
\begin{gather*}
X:\;
\left(
\begin{array}{ c c c c}
\theta_0 & 0 & 0 & 0   \\
1 & \theta_1 & 0 & 0  \\
0 & 0 & \theta_1  & 0 \\
0 & 1 & 0  & \theta_2
\end{array}
\right),
\qquad
Y:\;
\left(
\begin{array}{ c c c c}
\theta_0 & \vartheta & q & 0   \\
0 & \theta_1 & 0 & 0  \\
0 & 0 & \theta_1  & 1 \\
0 & 0 & 0  & \theta_2
\end{array}
\right).
\end{gather*}
Consider the matrix that represents
$\xi_1X\xi_2-\xi_2 X \xi_1 $ with respect to the above basis.
For this matrix the $(4,3)$-entry is $-q^2$.
This entry is nonzero
so
$\xi_1 X \xi_2- \xi_2 X \xi_1$ is nonzero.
Interchanging the roles of $X$, $Y$ in the above argument,
we see that
 $\xi_1 Y \xi_2- \xi_2 Y \xi_1$ is
nonzero.
The result follows.
\end{proof}

\begin{Theorem}
The homomorphism $\phi: {\mathcal O}\to \Delta$ from
Definition~{\rm \ref{def:hatvp}}
is not injective.
\end{Theorem}

\begin{proof}
Combine
Lemma~\ref{def:hatvpbasic}
and
Lemma~\ref{lem:nonz}.
\end{proof}

\section[The 2-sided ideal $\Delta\lbrack  \Delta,\Delta\rbrack \Delta $]{The 2-sided ideal $\boldsymbol{\Delta\lbrack  \Delta,\Delta\rbrack \Delta}$}

We will be discussing the following
subspace of $\Delta $:
\begin{gather*}
\lbrack \Delta ,\Delta \rbrack =
{\rm Span}\lbrace \lbrack u,v \rbrack \,|\,u,v \in \Delta \rbrace.
\end{gather*}
Observe that
$\Delta \lbrack  \Delta ,\Delta \rbrack \Delta $
is the 2-sided ideal of $\Delta $ generated by
$\lbrack \Delta ,\Delta \rbrack$.
In this section we describe this ideal from several
points of view.

 Let
$\overline A$,
$\overline B$,
$\overline C$ denote mutually commuting indeterminates.
Let $\mathbb F\lbrack
\overline A,
\overline B,
\overline C
\rbrack$
denote the $\mathbb F$-algebra consisting of the polynomials
in $
\overline A$,
$\overline B$,
$\overline C
$
 that have all coef\/f\/icients in $\mathbb F$.

\begin{Lemma}
\label{lem:com2}
There exists a unique $\mathbb F$-algebra homomorphism
$\Delta \to
\mathbb F\lbrack
\overline A,
\overline B,
\overline C
\rbrack $
that sends
\begin{gather*}
A \mapsto \overline A,
\qquad
B \mapsto \overline B,
\qquad
C \mapsto \overline C.
\end{gather*}
This homomorphism is surjective.
\end{Lemma}

\begin{proof}
The algebra
$
\mathbb F\lbrack
\overline A,
\overline B,
\overline C
\rbrack $ is commutative, so each of
\begin{gather*}
 \overline A+ \frac{q\overline B \,
\overline C-q^{-1}\overline C \, \overline B}{q^2-q^{-2}},
\qquad
\overline B+
\frac{q\overline C\, \overline A-q^{-1}\overline A \,\overline C}{q^2-q^{-2}},
\qquad
\overline C+
\frac{q\overline A\,\overline B-q^{-1}\overline B\,\overline A}{q^2-q^{-2}}
\end{gather*}
is central in
$
\mathbb F\lbrack
\overline A,
\overline B,
\overline C
\rbrack $.
Consequently
$\overline A$,
$\overline B$,
$\overline C
$ satisfy the def\/ining relations for~$\Delta $ from
Def\/inition~\ref{def:uaw}. Therefore the homomorphism exists.
The homomorphism is unique since
$A$,
$B$,
$C$ generate~$\Delta $.
The homomorphism is surjective
since
$
\overline A$,
$\overline B$,
$\overline C
$
generate
$
\mathbb F\lbrack
\overline A,
\overline B,
\overline C
\rbrack $.
\end{proof}

\begin{Definition}
Referring to the map
$\Delta \to \mathbb F\lbrack
\overline A,
\overline B,
\overline C
\rbrack $
from Lemma~\ref{lem:com2}, for all
$u\in \Delta $ let $\overline u$ denote the image of
$u$.
\end{Definition}

\begin{Lemma}
\label{lem:com2a}
We have
\begin{alignat*}{3}
&\overline \alpha
=
\big(q+q^{-1}\big)\overline A
+
{\overline B}\, {\overline C},
\qquad &&
\overline \beta
=
\big(q+q^{-1}\big)\overline B
+
\overline C \,\overline A, &
\\
&
\overline \gamma
=
\big(q+q^{-1}\big)\overline C
+
\overline A \,\overline B,
\qquad &&  \overline \Omega =
-
\big(q+q^{-1}\big)\overline A \, \overline B \,\overline C
-
{\overline A}^2
-
{\overline B}^2
-
{\overline C}^2. &
\end{alignat*}
\end{Lemma}

\begin{proof}
The assertions about $\alpha$, $\beta$, $\gamma $
 follow from
Def\/inition~\ref{def:abc}.
The assertion about~$\Omega $ follows
from its def\/inition in
Lemma~\ref{lem:six}.
\end{proof}

\begin{Proposition}
\label{lem:coincide}
The following coincide:
\begin{enumerate}\itemsep=0pt
\item[$(i)$]
The $2$-sided ideal
$\Delta \lbrack  \Delta ,\Delta \rbrack \Delta $;
\item[$(ii)$]
The kernel of the homomorphism
$\Delta \to
\mathbb F\lbrack
\overline A,
\overline B,
\overline C
\rbrack $ from Lemma~{\rm \ref{lem:com2}}.
\end{enumerate}
\end{Proposition}

\begin{proof}
Let $\Gamma$ denote the kernel of the
 homomorphism
$\Delta \to
\mathbb F\lbrack
\overline A,
\overline B,
\overline C
\rbrack $
 and note that~$\Gamma $ is a~2-sided ideal
 of~$\Delta $.
We have
$\lbrack \Delta ,\Delta \rbrack \subseteq \Gamma$
since
$
\mathbb F\lbrack
\overline A,
\overline B,
\overline C
\rbrack $ is commutative,
and
$\Delta \lbrack \Delta ,\Delta \rbrack \Delta
\subseteq \Gamma $
since
$\Gamma $ is a 2-sided ideal of $\Delta $.
The elements
$\lbrace \overline A^i \overline B^j \overline C^k |i,j,k\geq 0\rbrace $
form a  basis for the $\mathbb F$-vector space
$
\mathbb F\lbrack
\overline A,
\overline B,
\overline C
\rbrack $.
Therefore
the elements $\lbrace A^iB^jC^k |i,j,k\geq 0\rbrace$
form a basis for a complement of
$\Gamma $ in $\Delta $.
Denote this complement by~$M$, so
the sum
$\Delta =M+\Gamma $ is direct.
We now show $\Delta =M+
\Delta \lbrack \Delta ,\Delta \rbrack \Delta $.
The $\mathbb F$-algebra~$\Delta $
is generated by~$A$,~$B$,~$C$. Therefore
the $\mathbb F$-vector space~$\Delta $
is spanned by elements of the form
$x_1x_2\cdots x_n$ where $n\geq 0$ and
$x_i \in \lbrace A,B,C\rbrace$ for
$1 \leq i \leq n$.
Let
$x_1x_2\cdots x_n$ denote such
an element. Then
$\Delta \lbrack \Delta ,\Delta \rbrack \Delta $ contains
\begin{gather*}
x_1\cdots x_{i-1}x_i\cdots x_n  -
x_1\cdots x_ix_{i-1}\cdots x_n
\end{gather*}
for $2 \leq i \leq n$.
Therefore
$\Delta \lbrack \Delta ,\Delta \rbrack \Delta $ contains
\begin{gather*}
x_1x_2\cdots x_n  -  A^iB^j C^k,
\end{gather*}
where $i$, $j$, $k$ denote the number of
times $A$, $B$, $C$ appear among
$x_1, x_2,\ldots,  x_n$.
Observe $A^iB^jC^k \in M$ so
$x_1x_2\cdots x_n \in
M+
\Delta \lbrack \Delta ,\Delta \rbrack \Delta $.
Therefore
$\Delta =M+
\Delta \lbrack \Delta ,\Delta \rbrack \Delta $.
We already showed
$\Delta \lbrack \Delta ,\Delta \rbrack \Delta
\subseteq \Gamma $
and the sum
$\Delta =M+\Gamma $ is direct. By these comments
$\Delta \lbrack \Delta ,\Delta \rbrack \Delta
=\Gamma $.
\end{proof}

Recall
the ${\rm  {PSL}}_2(\mathbb Z)$-action on $\Delta $ from Theorem
\ref{thm:try}.
We now relate
this action
to the
homomorphism
$\Delta \to
\mathbb F\lbrack \overline A,\overline B, \overline C\rbrack$
from Lemma~\ref{lem:com2}.

\begin{Lemma}
\label{lem:overlineaction}
The group
${\rm  {PSL}}_2(\mathbb Z)$ acts on
$\mathbb F\lbrack \overline A,\overline B, \overline C\rbrack$
as a group of automorphisms in the following way:
\[
\begin{array}{@{}c| ccc }
u &  \overline A & \overline B  & \overline C
\\
\hline
\vphantom{\big|}\rho(u) &  \overline B & \overline C & \overline A
\\
\sigma(u)  &  \overline B & \overline A & \overline C
\end{array}
\]
\end{Lemma}

\begin{proof}
$\mathbb F\lbrack \overline A,\overline B, \overline C\rbrack$
has an automorphism of order 3 that sends
$(\overline A, \overline B, \overline C)\to
(\overline B, \overline C, \overline A) $, and an automorphism
of order 2 that sends
$(\overline A, \overline B, \overline C)\to
(\overline B, \overline A, \overline C)$.
\end{proof}

\begin{Lemma}
\label{lem:comdiagoverline}
For $g \in
{\rm  {PSL}}_2(\mathbb Z)$
 the following diagram commutes:
\[
\begin{CD}
\Delta  @>u \mapsto \overline u>>
         \mathbb F\lbrack \overline A, \overline B, \overline C\rbrack  \\
	  @Vg VV                     @VVg V \\
	 \Delta
	  @>>u \mapsto \overline u>
         \mathbb F\lbrack \overline A, \overline B,\overline C\rbrack
		    \end{CD}
		    \]
\end{Lemma}

\begin{proof}
Without loss
we may assume that
  $g$ is one of $\rho$, $\sigma$.
By Theorem~\ref{thm:try}
the action of $\rho$ on $\Delta $
cyclically permutes
$A$, $B$, $C$.
By Lemma~\ref{lem:overlineaction}
the action of $\rho$ on
$\mathbb F\lbrack \overline A,\overline B,\overline C\rbrack $
cyclically permutes
$\overline A$, $\overline B$, $\overline C$.
By Theorem~\ref{thm:try}
the action of~$\sigma$ on~$\Delta$
swaps $A,B$ and f\/ixes~$\gamma$. By Lemma~\ref{lem:com2a} and
Lemma~\ref{lem:overlineaction},
the action of $\sigma $ on
$\mathbb F\lbrack \overline A,\overline B,\overline C\rbrack $
swaps $\overline A$, $\overline B$
and f\/ixes~$ \overline \gamma$.
By these comments the diagram commutes.
\end{proof}

\begin{Definition}
\label{def:P}
By Lemma~\ref{lem:overlineaction}
each element of
${\rm  {PSL}}_2(\mathbb Z)$
permutes $\overline A$, $\overline B$, $\overline C$. This
induces a~group homomorphism
from ${\rm  {PSL}}_2(\mathbb Z)$ onto the symmetric group~$S_3$.
Let $\mathbb P$ denote the  kernel of this homomorphism.
Thus $\mathbb P$ is a normal subgroup of
${\rm  {PSL}}_2(\mathbb Z)$
and the quotient group
${\rm  {PSL}}_2(\mathbb Z)/\mathbb P$ is isomorphic to~$S_3$.
\end{Definition}

Our last main goal is to show that
\begin{gather}
\label{eq:fnlgl}
\Delta \lbrack \Delta ,\Delta \rbrack \Delta + \mathbb F 1 =
\langle A,B\rangle
\cap
\langle B,C\rangle
\cap
\langle A,C\rangle.
\end{gather}
 Note that in~(\ref{eq:fnlgl}) the sum on the left is
direct; otherwise
the ideal
$\Delta \lbrack \Delta ,\Delta \rbrack \Delta$ contains 1
and is therefore equal to $\Delta$, contradicting
Proposition
\ref{lem:coincide}.

\begin{Definition}
\label{def:hatu}
For notational convenience abbreviate
$\mathbb O =\langle A,B\rangle$. Thus
$\mathbb O$ is the image of~$\mathcal O$ under
the homomoprhism $\phi$
from Def\/inition~\ref{def:hatvp}.
\end{Definition}

\begin{Lemma}
\label{lem:notdir}
$\Delta = \sum_{n=0}^\infty \mathbb O \gamma^n$.
\end{Lemma}

\begin{proof}
 The
 algebra $\Delta $ is generated by
 $\mathbb O$, $\gamma$.
Moreover~$\gamma$ is central in~$\Delta $.
\end{proof}

\begin{Note}  The sum in Lemma
\ref{lem:notdir}
is not direct, by the third displayed equation in
Theorem~\ref{prop:altpres}.
\end{Note}

Note that
$\mathbb O\lbrack  A,B\rbrack \mathbb  O$
is the 2-sided ideal of
$\mathbb  O$ generated by
$\lbrack  A,B\rbrack$.

\begin{Lemma}
\label{lem:three}
The following $(i)$--$(iii)$ hold.
\begin{enumerate}\itemsep=0pt
\item[$(i)$]
$\lbrack \mathbb O,\mathbb O\rbrack \subseteq
\mathbb O \lbrack A,B\rbrack \mathbb O$.
\item[$(ii)$]
$  \lbrack A,B\rbrack
\gamma \in
\lbrack \mathbb O,\mathbb O\rbrack $.
\item[$(iii)$]
$\mathbb O \lbrack A,B\rbrack \mathbb O \gamma \subseteq
\mathbb O \lbrack A,B\rbrack \mathbb O$.
\end{enumerate}
\end{Lemma}
\begin{proof}\sloppy
$(i)$
Abbreviate $R=
\mathbb O \lbrack A,B\rbrack \mathbb O$ and consider the
quotient algebra $\mathbb O/R$.
The ele\-ments~$A$,~$B$ generate~$\mathbb O$, and
these generators satisfy
$\lbrack A,B\rbrack \in R$.
Therefore $A+R$, $B+R$ generate
$\mathbb O/R$, and these generators commute.
This shows that $\mathbb O/R$ is commutative.
Consider the canonical map $\mathbb O \to \mathbb O/R$.
This map has kernel~$R$. The map sends
$\lbrack \mathbb O,\mathbb O\rbrack \mapsto 0 $ since~$\mathbb O/R$ is commutative. Therefore
$\lbrack \mathbb O,\mathbb O\rbrack \subseteq R$.

  $(ii)$ In the third displayed equation of
Theorem~\ref{prop:altpres},
the expression on the right is a nonzero
scalar multiple of
$\lbrack A,B \rbrack \gamma $.
The expression on the left
is equal to
$\lbrack A^2,B^2 \rbrack +(q^2+q^{-2})
\lbrack B,ABA \rbrack $
and is therefore in
$\lbrack \mathbb O,\mathbb O\rbrack $.
The result follows.

  $(iii)$ By $(i)$, $(ii)$ above and
since $\gamma$ is central.
\end{proof}

\begin{Lemma}
\label{lem:ideal}
$\mathbb O \lbrack A,B\rbrack \mathbb O$ is a $2$-sided ideal of
$\Delta $.
\end{Lemma}

\begin{proof}
Abbreviate $R =
\mathbb O \lbrack A,B\rbrack \mathbb O$.
We show $\Delta R\subseteq R$ and
$R\Delta \subseteq R$.
Recall that $\Delta $ is generated by
$\mathbb O$, $\gamma$.
By construction
 $\mathbb O R\subseteq R$ and
$R \mathbb  O\subseteq R$.
By Lemma~\ref{lem:three}$(iii)$ and
since $\gamma$ is central we have
 $\gamma R\subseteq R$ and
$R \gamma \subseteq R$.
By these comments
$\Delta R\subseteq R$ and
$R\Delta \subseteq R$.
\end{proof}

\begin{Lemma}
\label{lem:comu}
We have
\begin{gather*}
\mathbb O \lbrack A,B \rbrack \mathbb O=
\Delta\lbrack  \Delta,\Delta\rbrack \Delta .
\end{gather*}
\end{Lemma}

\begin{proof}
We  have
$\mathbb O \subseteq \Delta $
and
$\lbrack A,B\rbrack \in \lbrack \Delta ,\Delta \rbrack$
so
$\mathbb O \lbrack A,B \rbrack \mathbb  O\subseteq
\Delta\lbrack  \Delta,\Delta\rbrack \Delta $.
We now show the reverse inclusion.
To this end we analyze
$\lbrack \Delta ,\Delta \rbrack$
using
Lemma~\ref{lem:notdir}.
For integers $m,n\geq 0$,
\begin{gather*}
\lbrack \mathbb  O \gamma^m,\mathbb O \gamma^n\rbrack
 =
\lbrack \mathbb  O ,\mathbb  O \rbrack \gamma^{m+n}
\overset{\text{by Lemma \ref{lem:three}({\it i})}}{\subseteq}
\mathbb  O \lbrack A,B\rbrack \mathbb O \gamma^{m+n}
 \overset{\text{by Lemma \ref{lem:ideal}}}{\subseteq}
\mathbb  O \lbrack A,B\rbrack \mathbb O.
\end{gather*}
By
this and Lemma~\ref{lem:notdir} we
obtain
$\lbrack \Delta ,\Delta \rbrack \subseteq
\mathbb  O \lbrack A,B\rbrack \mathbb O$.
Now using
Lemma~\ref{lem:ideal} we obtain
$\Delta \lbrack \Delta ,\Delta \rbrack \Delta  \subseteq
\mathbb  O \lbrack A,B\rbrack \mathbb  O$.
The result follows.
\end{proof}

In the algebra $\mathbb F\lbrack \overline A,\overline B, \overline C\rbrack$
let
$\mathbb F\lbrack \overline A,\overline B\rbrack$
(resp.~$\mathbb F\lbrack \overline B,\overline C\rbrack$)
(resp.~$\mathbb F\lbrack \overline A,\overline C\rbrack$)
denote the subalgebra generated by~$\overline A$,~$\overline B$
(resp.~$\overline B$, $\overline C$)
(resp.~$\overline A$, $\overline C$).

\begin{Proposition}
\label{thm:nice}
Referring to the homomorphism
$\Delta \to
\mathbb F\lbrack \overline A,\overline B, \overline C\rbrack$
from Lemma~{\rm \ref{lem:com2}},
\begin{enumerate}\itemsep=0pt
\item[$(i)$]
$\langle A,B \rangle $ is the preimage of
$\mathbb F\lbrack \overline A,\overline B\rbrack$;
\item[$(ii)$]
$\langle B,C\rangle $ is the  preimage of
$\mathbb F\lbrack \overline B,\overline C\rbrack$;
\item[$(iii)$]
$\langle A,C\rangle $ is the  preimage of
$\mathbb F\lbrack \overline A,\overline C\rbrack$.
\end{enumerate}
\end{Proposition}

\begin{proof}
$(i)$ Recall $\mathbb O = \langle A,B\rangle $.
For the homomorphism
in the proposition statement
the image of
$\mathbb O$ is
$\mathbb F\lbrack \overline A,\overline B\rbrack$.
Therefore the preimage of
$\mathbb F\lbrack \overline A,\overline B\rbrack$ is
$\mathbb O$ plus the kernel.
The kernel is
$\Delta\lbrack  \Delta,\Delta\rbrack \Delta $
by Proposition~\ref{lem:coincide}, and
this is contained in
$\mathbb O$ by
Lemma~\ref{lem:comu}.
Therefore~$\mathbb O $ is the preimage of
$\mathbb F\lbrack \overline A,\overline B\rbrack$.

  $(ii)$, $(iii)$
Apply $\rho $ twice
 to everything
in part $(i)$ above.
\end{proof}

\begin{Theorem}
\label{thm:big}
$\Delta\lbrack  \Delta,\Delta\rbrack \Delta + \mathbb F 1 =
\langle A,B\rangle
\cap
\langle B,C\rangle
\cap
\langle A,C\rangle$.
\end{Theorem}

\begin{proof}
By
Proposition~\ref{lem:coincide},
Proposition~\ref{thm:nice}, and since
$\mathbb F\lbrack \overline A,\overline B\rbrack
\cap
\mathbb F\lbrack \overline B,\overline C\rbrack
\cap
\mathbb F\lbrack \overline A,\overline C\rbrack = \mathbb F 1$.
\end{proof}

We f\/inish with some comments related to
Proposition~\ref{thm:nice}
and Theorem~\ref{thm:big}.

\begin{Proposition}\label{prop:complement1} \sloppy
In the table below, each space~$U$ is a
 subalgebra of $\Delta $ that
contains
$\Delta\lbrack  \Delta,\Delta\rbrack \Delta$.
The elements to the right of~$U$
form a basis for
a complement of
$\Delta\lbrack  \Delta,\Delta\rbrack \Delta$ in~$U$.
\begin{center}
\begin{tabular}{c| c}
$U$ &  basis for a complement of
$\Delta\lbrack  \Delta,\Delta\rbrack \Delta$ in $U$\bsep{1pt}
\\
\hline \hline
$\Delta$ &  $A^iB^jC^k\qquad i,j,k\geq 0$\tsep{1pt}
\\
\hline
$\langle A,B\rangle $ &  $A^iB^j\qquad i,j\geq 0$\tsep{1pt}
\\
$\langle B,C\rangle $ &  $B^jC^k\qquad j,k\geq 0$
\\
$\langle A,C\rangle $ &  $A^iC^k\qquad i,k\geq 0$
\\
\hline
$\langle A,B\rangle \cap
\langle A,C\rangle
$ &  $A^i\qquad i\geq 0$\tsep{1pt}
\\
$\langle A,B\rangle \cap
\langle B,C\rangle
$ &  $B^j\qquad j\geq 0$
\\
$\langle A,C\rangle \cap
\langle B,C\rangle
$ &  $C^k\qquad k\geq 0$
\\
\hline
$\langle A,B\rangle \cap
\langle B,C\rangle \cap
\langle A,C\rangle
$ &  $1\qquad \qquad $\tsep{1pt}
\end{tabular}
        \end{center}
\end{Proposition}

\begin{proof}
By Proposition
\ref{lem:coincide} and
Proposition \ref{thm:nice}.
\end{proof}

\begin{Proposition}
\label{prop:groupaction}
The automorphisms $\rho$ and $\sigma$
permute
the subalgebras
\begin{gather}
\label{eq:3set}
\langle A,B\rangle, \qquad
\langle B,C\rangle,\qquad
\langle A,C\rangle
\end{gather}
in the following way:
\[
\begin{array}{@{}c|ccc }
U &
\langle A,B\rangle  & \langle B,C\rangle   & \langle A,C\rangle
\\
\hline
\rho(U) &
\langle B,C\rangle  & \langle A,C\rangle   & \langle A,B\rangle
\\
\sigma(U)  &
\langle A,B\rangle  & \langle A,C\rangle   & \langle B,C\rangle
\end{array}
\]
\end{Proposition}

\begin{proof}
By Lemmas
\ref{lem:overlineaction},
\ref{lem:comdiagoverline}
and
Proposition \ref{thm:nice}.
\end{proof}

 By
Proposition
\ref{prop:groupaction}
the action of
${\rm  {PSL}}_2(\mathbb Z)$ on $\Delta$
induces an action of
${\rm  {PSL}}_2(\mathbb Z)$ on the 3-element
set~(\ref{eq:3set}).
The kernel of this action is
the
group $\mathbb P$ from
Def\/inition~\ref{def:P}.

\begin{Corollary}
\label{prop:P}
Each of the subalgebras
\begin{gather*}
\langle A,B\rangle, \qquad
\langle B,C\rangle,\qquad
\langle A,C\rangle
\end{gather*}
is invariant under the group $\mathbb P$ from Definition~{\rm \ref{def:P}}.
\end{Corollary}

\section{Directions for further research}

  In this section we give some suggestions
for further research.
Recall the algebra $\Delta $  from
Def\/inition~\ref{def:uaw}.

\begin{Problem}
Recall from the Introduction that $\Delta$
was originally motivated by the Askey--Wilson polynomials.
These polynomials are the most general
family
in a master class of ortho\-go\-nal polynomials called the
Askey scheme~\cite{KoeSwa}. For each polynomial
family
in the Askey scheme, there should be an analog of~$\Delta$
obtained from the appropriate version of~AW(3) by
interpreting parameters as central elements.
Investigate these other algebras along
the lines of the present paper.
\end{Problem}

\begin{Problem}
For this problem assume the characteristic of $\mathbb F$ is not~2.
By
Theorem~\ref{thm:try}
and Theorem~\ref{thm:faith} the group
${\rm  {PSL}}_2(\mathbb Z)$ acts faithfully on
$\Delta $ as a group of
automorphisms. This action induces an
injection of groups
${\rm  {PSL}}_2(\mathbb Z) \to {\rm Aut}(\Delta)$.
This injection is not an isomorphism
for the following reason.
Given any element of $\Delta$ among
$A$, $B$, $C$
there exists a unique automorphism of~$\Delta$
that f\/ixes that element and changes the sign of
the other two elements.
This automorphism is not contained in
the image of the above injection, because
its induced action on
$\mathbb F \lbrack  \overline A,\overline B,\overline C\rbrack$
does not match the description given in
Def\/inition~\ref{def:P}.
The above three automorphisms
are the nonidentity elements in
a subgroup $\mathbb K
\subseteq {\rm Aut}(\Delta)$ that is isomorphic
to the Klein 4-group $\mathbb Z_2 \times \mathbb Z_2$.
Do $\mathbb K$ and
${\rm  {PSL}}_2(\mathbb Z)$ together generate
${\rm Aut}(\Delta)$?
\end{Problem}

\begin{Problem}
View the $\mathbb F$-vector space $\Delta$ as a
${\rm  {PSL}}_2(\mathbb Z)$-module.  Describe the
irreducible
${\rm  {PSL}}_2(\mathbb Z)$-submodules of $\Delta$.
Is $\Delta $ a direct sum of irreducible
${\rm  {PSL}}_2(\mathbb Z)$-submodules?
\end{Problem}

\begin{Problem}
Find all the 2-sided ideals of $\Delta$.
Which of these are
${\rm  {PSL}}_2(\mathbb Z)$-invariant?
\end{Problem}

\begin{Problem}
Find all the
${\rm  {PSL}}_2(\mathbb Z)$-invariant subalgebras of
$\Delta$.
\end{Problem}

\begin{Problem}
\label{prob:pp}
Describe the subalgebra of
$\Delta $ consisting of the elements
in~$\Delta$ that are f\/ixed by everything in
the group~$\mathbb P$ from
Def\/inition~\ref{def:P}.
This subalgebra contains $\langle \Omega,\alpha,\beta,\gamma\rangle$.
Is this containment proper?
\end{Problem}

\begin{Problem}
Describe the
${\rm  {PSL}}_2(\mathbb Z)$-submodule of $\Delta$
that is generated by $\langle A \rangle $.
Also,
describe the $\mathbb P$-submodule of
 $\Delta$
that is generated by $\langle A \rangle $.
\end{Problem}

\begin{Problem}
Consider the basis for $\Delta$ given in
Theorem \ref{thm:basis} or
Theorem \ref{lem:basis2}.
Find the matrices that represent $\rho$ and $\sigma$
with respect to this basis.
Find the matrices that represent left-multiplication
by $A$, $B$, $C$ with respect to this basis.
Hopefully the entries in the above matrices are attractive
in some way.
If not, then f\/ind a  basis for~$\Delta $
with respect to which the above matrix entries are
attractive.
\end{Problem}

\begin{Problem}
Find a basis for the center~$Z(\Delta)$ under the assumption~$q$ is a root of unity.
\end{Problem}

\begin{Problem}
Give a basis for the $\mathbb F$-vector space
$\Delta \lbrack \Delta, \Delta \rbrack \Delta $.
\end{Problem}

\begin{Problem}
Recall the homomorhism $\Delta \to \mathbb F\lbrack \overline A,\overline B,
\overline C\rbrack $ from Lemma~\ref{lem:com2}. Restrict  this homomorphism
to $Z(\Delta)$ or $\langle A,Z(\Delta)\rangle$. In
each case
f\/ind a basis for the kernel and
image.
\end{Problem}

\begin{Problem}
Each of the following is a commutative subalgebra of~$\Delta$;
for each one give
a~basis and also a~presentation by generators and
relations.
\begin{enumerate}\itemsep=0pt
\item[$(i)$]
The intersection of
$\langle A,B\rangle$ and
$Z(\Delta)$.
\item[$(ii)$]
The intersection of $\langle A,B\rangle$ and
$\langle A, Z(\Delta) \rangle$.
\item[$(iii)$]
The intersection of
$\langle A,B\rangle$ and
$\langle C, Z(\Delta) \rangle$.
\end{enumerate}
\end{Problem}

\begin{Problem}
Find a Hopf algebra structure for $\Delta $.
See
\cite[Proposition~4.1]{basnc} and
\cite[Theorem~6.10]{ciccoli}
for some results in this direction.
\end{Problem}

Motivated by Theorem~\ref{thm:try}, let us view
$\Delta$ as a Lie algebra with
Lie bracket $\lbrack u,v\rbrack = uv-vu$ for
all $u,v \in \Delta$.

\begin{Problem}
Let $L$ denote the Lie subalgebra of $\Delta$
generated by $A$, $B$, $C$. Show that
\begin{gather*}
L \subseteq \mathbb F A + \mathbb F B + \mathbb F C +
\Delta\lbrack \Delta,\Delta\rbrack \Delta.
\end{gather*}
Show that $L$ is
${\rm  {PSL}}_2(\mathbb Z)$-invariant.
Find a basis for the
$\mathbb F$-vector space $L$. Describe
$L\cap Z(\Delta)$. Give a presentation for $L$
by generators and relations.
\end{Problem}

\pdfbookmark[1]{References}{ref}
\LastPageEnding

\end{document}